# ON RANDOM TOMOGRAPHY WITH UNOBSERVABLE PROJECTION ANGLES[1]

By Victor M. Panaretos

*Ecole Polytechnique Fédérale de Lausanne*

*To the memory of David George Kendall (1918–2007)
whose work will always be a source of inspiration*

We formulate and investigate a statistical inverse problem of a random tomographic nature, where a probability density function on $\mathbb{R}^3$ is to be recovered from observation of finitely many of its two-dimensional projections in random and unobservable directions. Such a problem is distinct from the classic problem of tomography where both the projections and the unit vectors normal to the projection plane are observable. The problem arises in single particle electron microscopy, a powerful method that biophysicists employ to learn the structure of biological macromolecules. Strictly speaking, the problem is unidentifiable and an appropriate reformulation is suggested hinging on ideas from Kendall's theory of shape. Within this setup, we demonstrate that a consistent solution to the problem may be derived, without attempting to estimate the unknown angles, if the density is assumed to admit a mixture representation.

## 1. Introduction.

The classical problem of tomography can be informally described as that of the determination of an object by knowledge of its projections in multiple directions. Problems of this nature arise in a variety of disciplines including medicine, astronomy, optics, geophysics and electron microscopy. Mathematically, the problem is formulated as that of seeking a solution to an integral equation relating a real function $f : \mathbb{R}^n \to \mathbb{R}$ to its

Received June 2008; revised November 2008.

[1]Supported in part by a NSF Graduate Research Fellowship and an NSF VIGRE Fellowship.

one-dimensional *Radon transform* (or *X-ray transform*),

$$(1.1) \qquad \check{f}(\xi, x) = \int_{-\infty}^{+\infty} f(x + \tau\xi) \, d\tau, \qquad \xi \in \mathbb{S}^{n-1} \text{ and } x \in \xi^{\perp}.$$

Under regularity conditions on $f$, the Radon transform can be seen to be invertible, and the function $f$ can be recovered by means of explicit formulas which we omit (Helgason [22] and Jensen [26]).

In practical situations, such as X-ray imaging (e.g., Shepp and Kruskal [47]), one seeks to determine $f$ given finitely many pairs of orientations and corresponding profiles $\{(\xi_i, \check{f}(\xi_i, \cdot))\}_{i=1}^{N}$. Several algorithms have been proposed to address this problem, and these are often problem specific, although one may single out broad classes, such as Fourier methods (based on the *projection-slice theorem*) and *back-projection* methods (see Natterer [38]). The subject matter and mathematical literature on such problems and their solution approaches is vast (see Deans [10] for a succinct overview).

In statistics, the tomographic reconstruction problem manifests itself most prominently in the case of *positron emission tomography* (*PET*), which can be seen as a special type of density estimation problem where a density $f$ is to be estimated given a random sample $\{(\Xi_i, X_i)\}_{i=1}^{n}$ from a density proportional to $\check{f}(\xi, x)$ (see Shepp and Vardi [48] and Vardi, Shepp and Kaufman [53]). PET lends itself to statistical treatment through a broad range of techniques such as likelihood-based, orthogonal series (singular value decomposition) and smoothed backprojection techniques, to name only a few (e.g., Vardi, Shepp and Kaufman [53], Silverman et al. [49], Green [20], Jones and Silverman [28]). Naturally, theoretical aspects such as consistency and optimality have also been widely investigated (e.g., Chang and Hsiung [8], Johnstone and Silverman [27] and O'Sullivan [39]). Further to PET, statistical problems such as random coefficient regression and multivariate density estimation have also been treated by means of insights and techniques gained from the field of tomography (Beran, Feuerverger and Hall [2], Feuerverger and Vardi [14] and O'Sullivan and Pawitan [40]).

In this paper, we formulate and investigate a stochastic variant of the classical tomographic reconstruction problem, where the profile orientations are not only random, but are in fact *unobservable*. This variant arises in the electron microscopy of single biological particles (see Section 1.1), a powerful method that biophysicists employ in order to study the structure of biological macrocolecules. It is qualitatively different from the usual tomography settings, where reconstruction crucially depends on the knowledge of the projection directions $\{\xi_i\}_{i=1}^{n}$. When the latter are unavailable, it is natural to ask whether anything interesting can be statistically said about the unknown density. We explore the limitations that are inherent when trying to answer such a question, and propose a mixture framework where the three-dimensional structure can be consistently estimated up to an orthogonal transformation, without attempting to estimate the unknown projection angles.

The paper is structured as follows. In Section 1.1 we present an informal introduction to the missing angle tomography problem of single particle electron microscopy. We then formulate the problem statistically (Section 2) and discuss its main aspects and the relevance of shape-theoretic ideas (Section 3). We then proceed to introduce a parametric framework (Section 3.1) which allows for a "statistical inversion" in the shape domain (Section 4). Illustrations are provided in Section 5 and the paper concludes with some remarks in Section 6.

1.1. *Single particle electron microscopy.* Resolving the structure of a biological particle is an undertaking that involves piecing together numerous facets of a complex investigation. The most important of these is, perhaps, the three-dimensional visualization of the particle whose dimension can be of the order of Angstroms (1 Å $= 10^{-10}$ m). Although it is X-rays that have traditionally been associated with particle structure determination, electron microscopy has arisen as a powerful tool with important advantages in these endeavors (Chiu [9], Frank [15], Glaeser et al. [19] and Henderson [23]).

The structure of a biological particle is described by the relative positioning of its atoms in space. Each atom's electrons create a potential around it, and the ensemble of these potentials gives rise to a potential distribution in three-dimensional space, the *shielded Coulomb potential distribution*, which is typically assumed to admit a *potential density* function, say $\rho(x, y, z)$. The structure of the particle is then described by this density.

This potential density provides the means of interaction with the electron microscope's beam: when the beam passes through the specimen (particle) in the $z$-direction, there is a reduction to the beam intensity caused by the scattering of electrons due to the interaction with the specimen. According to the *Abbe image formation theory* (Glaeser et al. [19]), the intensity recorded on the film under the specimen is approximately linear in the projection of the particle density in the $z$-direction, $\int_{-\infty}^{+\infty} \rho(x, y, z) \, dz$. Therefore, the imaging mode of the electron microscope provides us with a sample profile from the Radon transform of the particle's potential density.

While as such, the problem should be amenable to the "traditional" tomographic reconstruction techniques, things in practice are not as straightforward due to the problem of *radiation damage* (Glaeser [17]). Extended exposure to the electron beam will cause chemical bonds to break, and thus will alter the structure of the specimen. It follows that it is impossible to image the same particle under many different orientations. The exposure should instead be distributed over many identical particles. This can be achieved by crystallizing multiple particles (Drenth [12]) but reliance on crystals has several fundamental drawbacks (Frank [15] and Glaeser [18]). *Single particle cryo-electron microscopy* is a technique of electron microscopy, that aims at obtaining a three-dimensional representation of the particle without crystallizing the sample (e.g., Glaeser [18]). In this approach, a large number of structurally identical particles are imbedded unconstrained (i.e., uncrystallized) in an aqueous solution, then rapidly frozen and finally imaged via the electron microscope. Since the particles are unconstrained, they move and rotate freely within the aqueous solution, assuming haphazard orientations at the moment they are frozen. After preliminary processing, the data yielded are essentially noisy versions of the projected potential densities, at *random* and *unknown orientations*. Figures 1 and 2 present characteristic examples of such data in the presence of noise (Figure 1), and in the ideal—but practically impossible—noiseless case (Figure 2), for two different particles.

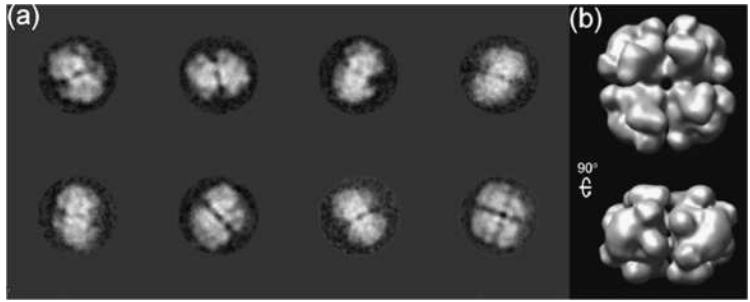

FIG. 1. (a) *Random profiles of pyruvate-ferredoxin oxidoreductase (PFOR) obtained via single-particle electron microscopy at the Lawrence Berkeley National Laboratory (Courtesy of F. Garczarek and R. M. Glaeser).* (b) *Reconstruction of the three-dimensional potential density after the projection angles have been estimated (Garczarek et al. [16]).*

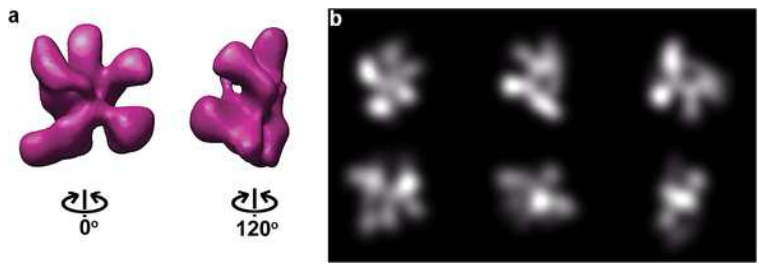

FIG. 2. (a) *A model of the three-dimensional potential density of the human translation initiation factor eIF3 derived from single particle data after angle estimation (Siridechadilok et al. [50]).* (b) *Noiseless random projections obtained from the known model (Courtesy of R. J. Hall).*

Biophysicists typically proceed via attempting to estimate the unobservable orientations, in order to then be able to iteratively use standard tomographic techniques (Frank [15] and Glaeser et al. [18]). However, they often rely on a priori knowledge on the structure of the particle either from other experiments or from an ad hoc examination of the projections by eye, in order to perform this estimation. Once an initial model is provided, it may be used to estimate the unknown angles and update the estimate. In cases where previous structural information is not available, and a naked eye examination is either not feasible (e.g., when the particle has no symmetries) or would best be avoided, it is natural to wonder whether an "objective" initial model can be extracted directly from the data, without attempting to estimate the unknown angles.

## 2. A stochastic Radon transform.

We may distinguish three important aspects in the random tomography problem that arises in single particle electron microscopy: (I) the samples of the Radon transform are obtained at *haphazard* orientations $\xi$ which are thought as random, (II) the physics of the data collection process allows for the possibility of within-plane rotations in a projection, (III) one *does not observe* the projection orientations.

These aspects become clear once we have a precise working definition and, for this reason, we define a random analogue to the Radon transform. Let $(\mathsf{SO}(3), \mathcal{D})$ be the measurable space of special orthogonal matrices, with $\mathcal{D}$ the Borel $\sigma$-algebra generated by the standard Riemannian metric on $\mathsf{SO}(3)$. Also, let $(L^2(\Delta_2), \mathcal{B})$ be the measurable space of square integrable functions on the disc $\Delta_2 := \{x \in \mathbb{R}^2 : \|x\| \leq \pi\}$, equipped with the Borel $\sigma$-algebra $\mathcal{B}$ generated by the $L^2$-norm.

Let $\rho : \mathbb{R}^3 \to [0, \infty)$ be probability density function centered at zero. Since the object of any tomographic probe is necessarily finite, we shall restrict our attention to densities that are supported on the ball $\Delta_3 := \{x \in \mathbb{R}^3 : \|x\| \leq \pi\}$ and essentially bounded (i.e., $\operatorname{ess\,sup} \rho < \infty$).

We define the *projection operator* of $\rho$ as the mapping $\Pi\{\rho\} : \mathsf{SO}(3) \to L^2(\Delta_2)$ given by

$$(2.1) \qquad (\Pi\{\rho\}(A))(x, y) := \int_{-\infty}^{+\infty} A\rho(x, y, z) \, dz \qquad \forall A \in \mathsf{SO}(3),$$

where $A\rho(\mathbf{x}) := \rho(A^{-1}\mathbf{x})$ for $\mathbf{x} \in \mathbb{R}^3$ and $A \in \mathsf{SO}(3)$. Given an element $A_0 \in \mathsf{SO}(3)$, the function $\Pi\{\rho\}(A_0)$ is the *projection* (or *profile*) of $\rho$ at orientation $A_0$. In particular, $\Pi\{\rho\}$ is well defined as a random element of $L^2(\Delta_2)$ if we equip $(\mathsf{SO}(3), \mathcal{D})$ with a probability measure.

PROPOSITION 2.1 (Measurability). *Let $\rho : \Delta_3 \to [0, \infty)$ be an essentially bounded probability density function centered at the origin. The projection operator $\Pi\{\rho\}$ is a measurable mapping from $(\mathsf{SO}(3), \mathcal{D})$ to $(L^2(\Delta_2), \mathcal{B})$.*

PROOF. Let $\theta : \Delta_3 \to \mathbb{R}$ be a continuous function, and let $A_n \overset{n \to \infty}{\Longrightarrow} A$ be a convergent sequence in $\mathsf{SO}(3)$. By continuity, it holds that $A_n \theta \to A\theta$ pointwise [recall that $A_n \theta(\mathbf{x}) = \theta(A_n^{-1}\mathbf{x})$]. Combining this fact with the bounded convergence theorem shows that $\lim_{n \to \infty}(\Pi\{\theta\}(A_n))(x, y) = (\Pi\{\theta\}(A))(x, y)$, for all $(x, y) \in \mathbb{R}^2$. The bounded convergence theorem then implies $\|\Pi\{\theta\}(A_n) - \Pi\{\theta\}(A)\|_2 \to 0$.

Now let $\rho$ be as in the assumptions of the theorem and let $\epsilon > 0$. By Lusin's theorem, there exists a continuous function $\theta_\epsilon : \Delta_3 \to \mathbb{R}$ such that $\mathsf{Leb}\{\mathbf{x} \in \Delta_3 : \rho(\mathbf{x}) \neq \theta_\epsilon(\mathbf{x})\} < \epsilon$. By the triangle inequality,

$$
\begin{aligned}
\|\Pi(\rho)(A_n) - \Pi(\rho)(A)\|_2 \leq\ & \|\Pi(\rho)(A_n) - \Pi(\theta_\epsilon)(A_n)\|_2 \\
& + \|\Pi(\theta_\epsilon)(A_n) - \Pi(\theta_\epsilon)(A)\|_2 \\
& + \|\Pi(\theta_\epsilon)(A) - \Pi(\rho)(A)\|_2.
\end{aligned}
$$

If we let $n \to \infty$, our earlier analysis shows that the second term will vanish. The first and third term on the right-hand side are bounded above by $\epsilon \cdot C$, for some finite $C \geq 0$ since $\rho$ is nonnegative and essentially bounded, while $\theta_\epsilon$ is bounded. Since the choice of $\epsilon$ is arbitrary, this establishes the continuity of $\Pi\{\rho\}$ with respect to the relevant topologies and its measurability with respect to the corresponding Borel $\sigma$-algebras. $\quad\square$

For $N \geq 1$, let $\{A_n\}_{n=1}^N$ be independent and identically distributed random elements of the special orthogonal group $\mathsf{SO}(3)$ distributed according to normalized Haar measure. We define the *stochastic Radon transform* of length $N$ of $\rho$ as the i.i.d. collection of random projections $\{\Pi\{\rho\}(A_n)\}_{n=1}^N$, taking values in the *sample space* $(L^2(\Delta_2), \mathcal{B})$. These realizations of independent projections are not coupled with the corresponding orientations, that is, we observe $\Pi\{\rho\}(A_n)$ but not $A_n$. For this reason, we suppress the dependence on $A_n$ whenever this does not cause confusion, and write $\check{\rho}_n$ for $\Pi\{\rho\}(A_n)$. From the classical statistical perspective, we observe that any centered essentially bounded density $\rho$ on $\Delta_3$ induces a probability measure $P_\rho$ on the measurable space $(L^2(\Delta_2), \mathcal{B})$ via

$$
P_\rho[B] = \Psi\{A \in \mathsf{SO}(3) : \Pi\{\rho\}(A) \in B\}, \qquad B \in \mathcal{B},
$$

with $\Psi$ denoting normalized Haar measure on $(\mathsf{SO}(3), \mathcal{D})$. The stochastic Radon transform of length $N$ of $\rho$ is then an i.i.d. random sample from the distribution $P_\rho$ (a collection of $N$ independent random fields with law $P_\rho$).

**3. Invariance and shape.** We wish to consider the recovery of a density given its stochastic Radon transform, that is, to investigate the feasibility of a *statistical inversion* of the stochastic transform. When seen as an estimation problem, the recovery problem exhibits certain special group invariance properties; these are manifested both in the parameter space as well as in the sample space, as *unidentifiability* and *sufficiency*, respectively.

Focusing first on the parameter space, we recall that a parametric family of models (distributions) $\{P_\theta\}$ with parameter space $\Theta$ is *identifiable* if the mapping $\theta \mapsto P_\theta$ is a bijection. The probability model parameterized in our case is the distribution of a random profile $\Pi\{\theta\}$ with $\theta \in \mathcal{F}$, where $\mathcal{F}$ is the set of essentially bounded probability densities supported on $\Delta_3$ that are centered at the origin [understood as a subset of the metric space $L^2(\Delta_3)$]. However, the following parametrization is not well defined, in the sense that it is unidentifiable: if $B \in \mathsf{O}(3)$ so that $B^\top B = I$ and $\det(B) \in \{-1, 1\}$, we may put $Q = \mathrm{diag}\{1, 1, \det(B)\}$, and observe that

$$\Pi\{B\theta\} = \int_{-\infty}^{+\infty} AB\theta(x, y, z)\, dz = \int_{-\infty}^{+\infty} ABQ\theta(x, y, \det(B) \cdot z)\, dz \overset{d}{=} \Pi\{\theta\}$$

for any $A \sim \mathrm{Haar}[\mathsf{SO}(3)]$, by right invariance of Haar measure.

It follows that the probability law $P_\rho$ induced on $L^2(\Delta_2)$ by the parameter $\rho \in \mathcal{F}$ is the same as the law $P_{B\rho}$ for any $B \in \mathsf{O}(3)$. This suggests that ideally we could only recover the original function modulo $\mathsf{O}(3)$, which leads to the need for a parametrization of the model in terms of those characteristics of the functions of $\mathcal{F}$ that are invariant under orthogonal transformations. Formally, let $G(\mathcal{F}) = \{\gamma_A : A \in \mathsf{O}(3)\}$ be the *group of rotations and reflections* on the function class $\mathcal{F}$, with action

$$(3.1) \qquad (\gamma_A f)(u) = f(A^{-1}u), \qquad A \in \mathsf{O}(3),\ f \in \mathcal{F}.$$

Define the *shape* $[f]$ of a function $f \in \mathcal{F}$ as its orbit under the action of $G(\mathcal{F})$

$$(3.2) \qquad [f] = \{\gamma(f) : \gamma \in G(\mathcal{F})\}.$$

Consequently, we call the quotient space $\mathcal{F}/G(\mathcal{F})$ the *shape space* of $\mathcal{F}$, and we denote it by $\Sigma\mathcal{F}$. While we saw that we cannot recover "more than $[\rho]$" from the stochastic Radon transform of $\rho$, we prove next that shape *can* be potentially recovered given a sample from the stochastic Radon transform.

THEOREM 3.1 (Singular identifiability). *Let $\mathcal{F}$ be the set of probability densities supported on $\Delta_3$ that are centered at the origin and are essentially bounded. For $\theta \in \mathcal{F}$, let $P_{[\theta]}$ denote the probability distribution induced on the sample space $(L^2(\Delta_2), \mathcal{B})$ by $[\theta] \in \Sigma\mathcal{F}$ via the stochastic Radon transform $\Pi\{h\}$ of any $h \in [\theta]$. Then, for any two distinct elements $[f], [g] \in \Sigma\mathcal{F}$, the measures $P_{[f]}$ and $P_{[g]}$ are mutually singular.*

To prove this theorem, we will make use of the following result on Radon transforms (see Proposition 7.8 in Helgason [22]).

PROPOSITION 3.1. *Let $f : \mathbb{R}^d \to \mathbb{R}$ be function of compact support and $\Xi$ be an infinite subset of the unit sphere $\mathbb{S}^2$. Then $f$ is determined by the collection $\{\check{f}_\xi\}_{\xi \in \Xi}$, where*

$$\check{f}_\xi(x) = \int_{-\infty}^{+\infty} f(x + \tau\xi)\, d\tau, \qquad x \in \xi^\perp.$$

PROOF OF THEOREM 3.1. Assume that $[f], [g] \in \Sigma\mathcal{F}$ and $[f] \neq [g]$. Since $P_{[g]} = P_g$ and $P_{[g]} = P_g$, it suffices to show that $P_f \perp P_g$. Since the shapes $[f]$ and $[g]$ are distinct, we have

$$f \neq Bg \qquad \forall B \in \mathsf{SO}(3).$$

It follows that given any $\Gamma \in \mathsf{SO}(3)$, the set $\{A \in \mathsf{SO}(3) : \Pi\{f\}(A) = \Pi\{\Gamma g\}(A)\}$ has Haar measure zero,

$$(3.3) \qquad \Psi\{A \in \mathsf{SO}(3) : \Pi\{f\}(A) = \Pi\{\Gamma g\}(A)\} = 0 \qquad \forall \Gamma \in \mathsf{SO}(3).$$

For if this were not the case we could find an uncountably infinite set $\Xi \subseteq \mathbb{S}^2$ such that $\check{f}_\xi = (\Gamma g)_\xi$ for all $\xi \in \Xi$, where

$$\check{h}_\xi(x) = \int_{-\infty}^{+\infty} h(x + \tau\xi)\, d\tau, \qquad x \in \xi^\perp.$$

So, by means of Proposition 3.1, we would conclude $f = \Gamma g$, contradicting our assumption.

Consider an arbitrary coupling of $P_f$ and $P_g$, that is, let $(\mathsf{SO}(3), \mathcal{D}, \Psi)$ be as before, and let $h : (\mathsf{SO}(3), \mathcal{D}) \to (\mathsf{SO}(3), \mathcal{D})$ be any measurable function such that $\Psi\{A \in \mathsf{SO}(3) : \Pi\{g\}[h(A)] \in \cdot\} = P_g[\cdot]$.

Initially, we assume that $h(\cdot)$ is continuous. Since $\mathsf{SO}(3)$ acts transitively on itself, we may represent $h$ as

$$h(A) = A\Gamma_A, \qquad \Gamma_A \in \mathsf{SO}(3).$$

By continuity of $h$, it follows that $A \mapsto \Gamma_A$ is also continuous.

Now, let $\{\mathbf{A}_n\}_{n \geq 1}$ be a monotone sequence of partitions of $\mathsf{SO}(3)$ that become finer as $n$ increases. That is, $\mathbf{A}_n$ partitions $\mathsf{SO}(3)$ into $n$ disjoint sets $\{\mathcal{A}_n^i\}_{i=1}^n$ with the property that for every $j \in \{1, \ldots, n+1\}$ there exists an $i_j \in \{1, \ldots, n\}$ such that $\mathcal{A}_{n+1}^j \subseteq \mathcal{A}_n^{i_j}$. Define a sequence of "simple" measurable functions

$$h_n(A) = A\Gamma_i^n \qquad \text{on the set } \mathcal{A}_n^i,$$

where $\Gamma_n^i$ is defined as

$$\Gamma_n^i := \operatorname*{arg\,min}_{\Gamma \in \{\Gamma_A \,:\, A \in \overline{\mathcal{A}}_n^i\}} \left\{ \min_{A \in \overline{\mathcal{A}}_n^i} \|\Pi\{f\}(A) - \Pi\{g\}(A\Gamma)\|_2 \right\}$$

for $\overline{\mathcal{A}}_n^i$ the closure of $\mathcal{A}_n^i$. The above is well defined by compactness of $\overline{\mathcal{A}}_n^i$ and continuity of $A \mapsto \Gamma_A$. Now, since $h$ is continuous, we have

$$h_n \to h \qquad \forall A \in \mathsf{SO}(3).$$

Hence, by continuity of the projection mapping and by the dominated convergence theorem, we have, for all $A \in \mathsf{SO}(3)$,

$$(3.4) \qquad \|\Pi\{f\}(A) - \Pi\{g\}(h_n(A))\|_2 \uparrow \|\Pi\{f\}(A) - \Pi\{g\}(h(A))\|_2.$$

That the convergence is monotone follows from the definition of $h_n$: the partition sequence is monotone and for $A \in \overline{\mathcal{A}}_{n+1}^j \subseteq \overline{\mathcal{A}}_n^{i_j}$ we have that $\{\Gamma_A : A \in \overline{\mathcal{A}}_{n+1}^j\} \subseteq \{\Gamma_A : A \in \overline{\mathcal{A}}_n^{i_j}\}$. We now proceed to define the sets

$$\mathcal{K}_n := \{A \in \mathsf{SO}(3) : \|\Pi\{f\}(A) - \Pi\{g\}(h_n(A))\|_2 > 0\},$$

$$\mathcal{K} := \{A \in \mathsf{SO}(3) : \|\Pi\{f\}(A) - \Pi\{g\}(h(A))\|_2 > 0\}.$$

By definition of $h_n$,

$$\mathcal{K}_n \subseteq \mathcal{K}_{n+1} \qquad \forall n \geq 1.$$

Therefore, $\mathcal{K}_n \uparrow \bigcup_{n=1}^{\infty} \mathcal{K}_n$, where

$$\begin{aligned}
\bigcup_{n=1}^{\infty} \mathcal{K}_n &= \bigcup_{n=1}^{\infty} \{A \in \mathsf{SO}(3) : \|\Pi\{f\}(A) - \Pi\{g\}(h_n(A))\|_2 > 0\} \\
&= \{A \in \mathsf{SO}(3) | \exists n \geq 1 : \|\Pi\{f\}(A) - \Pi\{g\}(h_n(A))\|_2 > 0\} \\
&= \{A \in \mathsf{SO}(3) : \|\Pi\{f\}(A) - \Pi\{g\}(h(A))\|_2 > 0\} \\
&= \mathcal{K}.
\end{aligned}$$

The penultimate equality follows from the monotone convergence given in (3.4). Continuity of $\Psi$ from below leads us to the conclusion

$$\Psi[\mathcal{K}] = \Psi\left[\bigcup_{n=1}^{\infty} \mathcal{K}_n\right] = \lim_{n \to \infty} \Psi[\mathcal{K}_n].$$

On the other hand, by definition of $h_n$, we have that, for all $n \geq 1$,

$$\begin{aligned}
\Psi[\mathcal{K}_n] &= \Psi\{A \in \mathsf{SO}(3) : \|\Pi\{f\}(A) - \Pi\{g\}(h_n(A))\|_2 > 0\} \\
&= \Psi\left[\bigcup_{i=1}^{n} \{A \in \mathsf{SO}(3) : \|\Pi\{f\}(A) - \Pi\{g\}(h_n(A))\|_2 > 0\} \cap \mathcal{A}_n^i\right] \\
&= \Psi\left[\biguplus_{i=1}^{n} \{A \in \mathcal{A}_n^i : \|\Pi\{f\}(A) - \Pi\{g\}(A\Gamma_n^i)\|_2 > 0\}\right] \\
&= \sum_{i=1}^{n} \Psi\{A \in \mathcal{A}_n^i : \|\Pi\{f\}(A) - \Pi\{\Gamma_n^i g\}(A)\|_2 > 0\}
\end{aligned}$$

$$= \sum_{i=1}^{n} \Psi[\mathcal{A}_n^i] = 1$$

by appealing to the first part of our proof (3.3). In summary, we have shown that

$$\Psi[\mathcal{K}] = \Psi\{A \in \mathsf{SO}(3) : \Pi\{f\}(A) \neq \Pi\{g\}(h(A))\} = 1.$$

Now consider the case where $h$ is measurable, but not continuous. We recall that Lusin's theorem guarantees that for any $\delta > 0$ there exists a closed set $\mathcal{H}_\delta$ (and hence compact in our case) and a continuous function $h_\delta$ such that $h = h_\delta$ on $\mathcal{H}_\delta$ while $\Psi(\mathsf{SO}(3) \setminus \mathcal{H}_\delta) < \delta$. Therefore, for arbitrary measurable $h$, and given any $\delta > 0$,

$$
\begin{aligned}
\Psi[\mathcal{K}] &= \Psi[\mathcal{K} \cap \mathcal{H}_\delta] + \Psi[\mathcal{K} \cap \mathcal{H}_\delta^\mathsf{c}] \\
&= \Psi\{A \in \mathcal{H}_\delta : \Pi\{f\}(A) \neq \Pi\{g\}(h(A))\} \\
&\quad + \Psi\{A \in \mathsf{SO}(3) \setminus \mathcal{H}_\delta : \Pi\{f\}(A) \neq \Pi\{g\}(h(A))\} \\
&= \Psi\{A \in \mathcal{H}_\delta : \Pi\{f\}(A) \neq \Pi\{g\}(h_\delta(A))\} \\
&\quad + \Psi\{A \in \mathsf{SO}(3) \setminus \mathcal{H}_\delta : \Pi\{f\}(A) \neq \Pi\{g\}(h(A))\} \\
&\geq \Psi\{A \in \mathcal{H}_\delta : \Pi\{f\}(A) \neq \Pi\{g\}(h_\delta(A))\} \\
&= \Psi[\mathcal{H}_\delta] \\
&\geq 1 - \delta.
\end{aligned}
$$

The choice of $\delta$ being arbitrary, we conclude that the event $\{\Pi\{f\} \neq \Pi\{g\}\}$ has probability 1 for an arbitrary coupling. Strassen's theorem now implies that the total variation distance between $P_f$ and $P_g$ is 1, which completes the proof. □

It follows that, under isotropic projection orientations, the unknown density *is identifiable* up to an orthogonal transformation, regardless of whether or not we observe the projection angles [in fact, this remains true if Haar measure $\Psi$ is replaced by any measure $\Psi'$ on $\mathsf{SO}(3)$ with $\Psi' \ll \Psi$].

Shape is not just crucial as a notion in the context of the parameter space only. Under isotropic projection orientations, any orthogonal transformation of the two-dimensional projection data contains the same information on the function-valued parameter. Letting $G(L^2(\Delta_2))$ denote the group of rotations and reflections on $L^2(\Delta_2)$, we define the shape of an element $f$ in the sample space $L^2(\Delta_2)$ as $[f] = \{\alpha(f) : \alpha \in G(L^2(\Delta_2))\}$. We equip the corresponding shape space (quotient space) $M := L^2(\Delta_2)/G(L^2(\Delta_2))$ with the Borel $\sigma$-algebra $\mathcal{M}$ generated by the quotient topology. This turns the quotient space into a measurable space, and the quotient mapping into a measurable mapping, that is, *a statistic*. The next proposition establishes that the shape statistic is sufficient with respect to the original shape (see page 85 of Schervish [46] for the definition of abstract sufficiency).

THEOREM 3.2. *Let $\mathcal{F}$ be the set of probability densities supported on $\Delta_3$ that are centered at the origin and are essentially bounded. For $\theta \in \mathcal{F}$, let $P_{[\theta]}$ denote the probability distribution induced on the sample space $(L^2(\mathbb{R}^2), \mathcal{B})$ by $[\theta] \in \Sigma\mathcal{F}$ via the stochastic Radon transform $\Pi(h)$ of any $h \in [\theta]$. The mapping $\Pi(h) \mapsto [\Pi(h)]$ is a sufficient statistic for the parameter $[\theta]$ and a maximal invariant statistic with respect to the group $G(L^2(\Delta_2))$.*

Before we prove Theorem 3.2, we prove a lemma and recall two results from measure theory.

LEMMA 3.1. *Let $[\theta] \in \Sigma\mathcal{F}$ and $P_{[\theta]}$ be the law of $\Pi\{h\}$, induced by $h \in [\theta]$. Then, given any $B \in \mathcal{B}$, $\gamma \in G(L^2(\Delta_2))$, it holds that $P_{[\theta]}\{B\} = P_{[\theta]}\{\gamma B\}$.*

PROOF. There exists a $W \in \mathsf{O}(2)$ such that for $A \sim \mathrm{Haar}[\mathsf{SO}(3)]$

$$\gamma[\Pi\{\theta\}(A)(x,y)] \stackrel{d}{=} \int_{-\infty}^{+\infty} \begin{pmatrix} W & \mathbf{0} \\ \mathbf{0}^\top & \det(W) \end{pmatrix} A\theta(x,y,z)\,dz \stackrel{d}{=} \Pi\{\theta\}(A)(x,y),$$

the last equality following from the left invariance of Haar measure. □

The next two results can be found in Lemma 1.6 and Theorem 2.29 of Kallenberg [29].

LEMMA 3.2. *Let $(M,d)$ be a metric space with topology $\mathcal{T}$ and Borel $\sigma$-algebra $\mathcal{A}$. Then, for any $D \subset M$, the metric space $(D,d)$ has topology $\mathcal{T} \cap D$ and Borel $\sigma$-algebra $\mathcal{A} \cap D$.*

PROPOSITION 3.2. *Let $G$ be a locally compact second countable Hausdorff group that acts transitively and properly on a locally compact second countable Hausdorff space $S$. Then, there exists, uniquely up to renormalization, a $G$-invariant Radon measure $\mu \neq 0$ on $S$.*

PROOF OF THEOREM 3.2. Maximal invariance follows immediately from the definition of shape as the orbit under the group of orthogonal transformations. To prove sufficiency, we note that the space $(L^2(\Delta_2), \mathcal{B})$ is a standard Borel space since it is complete and separable in the metric induced by the $L^2$-norm. It follows that there exists a regular conditional distribution $\nu(B|[\theta], m) \colon \mathcal{B} \times \Sigma\mathcal{F} \times M \to [0,1]$ for $\Pi(\theta)$ given $[\Pi(\theta)]$,

$$(3.5) \qquad \nu(B|[\theta], m) := P_{[\theta]}\{B|[\Pi(\theta)] = m\}, \qquad B \in \mathcal{B},\ m \in M.$$

Therefore, sufficiency will follow if we can show that $\nu(B|[\theta], m)$ is functionally independent of $[\theta]$, that is, $\nu(B|w, m) = r(B, m)$, $\forall w \in \Sigma \mathcal{F}$. We begin by observing that $\nu(m|w, m) = 1$. Therefore, $\nu(\cdot|w, m)$ can be viewed as a probability measure on $(m, \mathcal{B} \cap m)$, where $\mathcal{B} \cap m := \{m \cap A : A \in \mathcal{B}\}$. Now let $\mathcal{T}$ be the natural topology of $L^2(\Delta_2)$, so that $\mathcal{B} \cap m = \sigma(\mathcal{T}) \cap m = \sigma(\mathcal{T} \cap m)$ is the Borel $\sigma$-algebra of subsets of $m$, generated by the subspace topology $\mathcal{T} \cap m$ (Lemma 3.2). But $(m, \mathcal{T} \cap m)$ is a locally compact second countable Hausdorff space. Hence, by Proposition 3.2, there exists a unique Radon measure (up to constant multiples) $r(B, m)$ on $(m, \mathcal{B} \cap m)$ such that $r(B, m) = r(\gamma B, m)$, for all $\gamma \in G(L^2(\Delta_2))$ and $B \in \mathcal{B} \cap m$. But Lemma 3.1 implies that $\nu(B|w, m) = \nu(\gamma B|w, m)$ for all $\gamma \in G(L^2(\Delta_2))$ and all $B \in \mathcal{B} \cap m$, and $\nu$ is a probability measure. Consequently, it must be that $\nu(B|w, m) = \lambda r(B, m)$ for some $\lambda > 0$, which completes the proof. $\square$

It follows that our analysis should concentrate on the concept of shape. On the one hand, it is the shape of the unknown density that we seek to estimate; on the other hand, we should base our estimate on the shape of the projections, that being a sufficient statistic.

A systematic mathematical study of *shape* in the case of finitely many labeled points in Euclidean space was initiated by Kendall [31]; his motivation was the question of existence of alignments in megalithic stone monuments (Kendall and Kendall [32]). In Kendall's approach, shape is the collection of those characteristics of a labeled point pattern that are invariant under rotation, translation and scaling. The shape spaces induced have a manifold structure, and their geometry depends both on the number of points and the dimension of the ambient space (Le and Kendall [34]). A closely related concept of shape with a different "representation theory" was independently proposed by Bookstein [3], who was interested in biological applications and morphometrics. In Kendall's terminology, our version of shape would be called "*unoriented shape-and-size*," to stress the fact that $\mathrm{O}(3)$ is quotiented out while scalings are not. Kendall and Le [33] provide a compendious review of statistical shape theory.

Though shape spaces of *finite point patterns* are well understood and widely used in applied work, there is apparently no practical formulation of the *shape of a function*. An active field of research focuses on practical parameterizations of the shape of closed curves on the plane and in space (e.g., Younes [54] and Small and Le [51]), the principle motivation being computer vision. Such ideas do not appear useful, though, when attempting to find connections between the shape of a function and the shape of its integral transform. For this reason, we will hinge on Euclidean shape-theoretic ideas that will enable us to establish such connections, via an appropriate parametrization (see Panaretos [41], Panaretos [42] and Panaretos [43]).

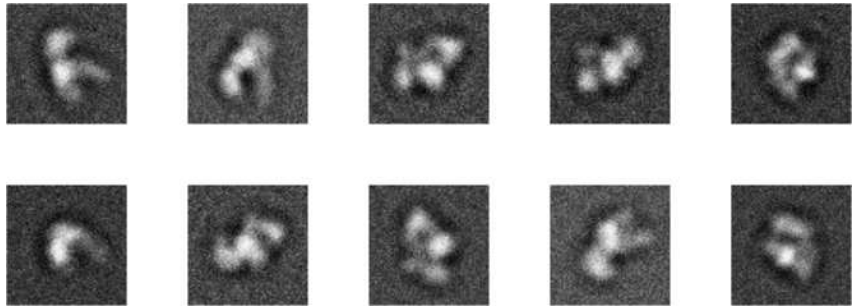

Fig. 3. *Synthetic single-particle projection data from the Klenow fragment of Escherichia coli DNA Polymerase I (Courtesy of A. Leschziner). The projections resemble mixtures of roughly circular components contaminated by noise.*

3.1. *Parametrization: radial expansions and the Gram matrix.* At least three ingredients come into play when considering a parametrization for the shape of a density in this particular context. First, it is important that the parametrization allow for the problem to be posed as one of parameter estimation. In addition, one may ask for a parametrization that makes it feasible to always explicitly be able to pick out a representative member from a particular shape class. Most important is the need to be able to find a connection between original shape and projected shape. In essence, what we are asking for is a parametrization that will allow us to convert the abstract setting of quotient spaces into something we can put a handle on.

With these general considerations in mind, we choose to focus on the following parametric yet flexible class of *finite mixtures of radial location densities*:

$$(3.6) \qquad \rho(\mathbf{x}) = \sum_{k=1}^{K} q_k \varphi(\mathbf{x}|\tilde{\mu}_k), \qquad \tilde{\mu}_k \in \mathbb{R}^3, \ q_k > 0, \ \sum_{k=1}^{K} q_k = 1,$$

where $\varphi(\cdot|\tilde{\xi})$ is a spherically symmetric probability density with expectation $\tilde{\xi}$, so that $\varphi(\cdot|\tilde{\xi}) = f(\|\mathbf{x} - \tilde{\xi}\|)$ for some probability density $f : \mathbb{R}^+ \to \mathbb{R}^+$. The choice of this particular type of expansion appears useful both from the applied and the mathematical points of view. The optics of the imaging procedure have a smoothing effect on the planar densities recorded on the film. As a result, the projected particles often do appear as an ensemble of roughly circular "blobs" (see, e.g., Figures 1, 2 and 3). Mixtures of Gaussians have previously been employed to obtain Riemannian metrics between biological shapes in deformation shape analysis (e.g., Peter and Rangarajan [44]).

From a mathematical point of view, this type of density is well behaved under orthogonal transformation and projection. For any $A \in \mathsf{O}(3)$,

$$(3.7) \qquad \varphi(A^\top \mathbf{x}|\tilde{\xi}) = f(\|A^\top \mathbf{x} - \tilde{\xi}\|) = f(\|\mathbf{x} - A\tilde{\xi}\|) = \varphi(\mathbf{x}|A\tilde{\xi}),$$

and letting $H$ be the projection onto the plane $z = 0$,

$$(3.8) \quad \int_{-\infty}^{+\infty} \varphi(x, y, z | A\tilde{\xi}) \, dz = \phi(x, y | \mu),$$

$$\begin{cases} \mu = HA\tilde{\xi}, \\ \phi(x, y | 0) = \int_{-\infty}^{+\infty} \varphi(x, y, z | 0) \, dz, \end{cases}$$

so that any rotation of the density can be encoded by a rotation of the location parameters $\tilde{\mu}_k$, while its two-dimensional profiles can be expressed as a mixture of the marginal of $\varphi$, regardless of the projection orientation.

To remove the effects of location, assume that the density is centered with respect to its location vectors, that is, assume $\sum_{k=1}^{K} \tilde{\mu}_k = 0$. Since any rotation of $\rho$ can be encoded by a rotation of its location vectors, we may use the characteristics of the location vectors to encode the shape of $\rho$. The *Gram matrix* generated by the collection $\{\tilde{\mu}_k\}$ is the $K \times K$ symmetric nonnegative definite matrix, whose $ij$th element is the inner product $\langle \tilde{\mu}_i, \tilde{\mu}_j \rangle$, as follows:

$$(3.9) \quad \mathsf{Gram}(\{\tilde{\mu}_k\}) = \begin{pmatrix} \|\tilde{\mu}_1\|^2 & \langle \tilde{\mu}_1, \tilde{\mu}_2 \rangle & \cdots & \langle \tilde{\mu}_1, \tilde{\mu}_K \rangle \\ \langle \tilde{\mu}_2, \tilde{\mu}_1 \rangle & \|\tilde{\mu}_2\|^2 & \cdots & \langle \tilde{\mu}_2, \tilde{\mu}_K \rangle \\ \vdots & & \ddots & \vdots \\ \langle \tilde{\mu}_K, \tilde{\mu}_1 \rangle & \cdots & & \|\tilde{\mu}_K\|^2 \end{pmatrix}.$$

In Kendall's shape theory, Gram matrices are employed as a coordinate system for the shape manifold induced by rigid motions. Note that if the vectors $\{\tilde{\mu}_k\}$ are arranged as the columns of a $3 \times K$ matrix $V$, then we may simply write $\mathsf{Gram}(V) = V^\top V$. The Gram matrix is invariant under orthogonal transformations of the generating vectors, since for $B \in \mathsf{O}(3)$ we immediately see that $\mathsf{Gram}(BV) = V^\top B^\top BV = V^\top V = \mathsf{Gram}(V)$. Furthermore, given a Gram matrix of rank $p$, one can find $K$ vectors in $\mathbb{R}^d$, $d \geq p$, with centroid zero whose pairwise inner products are given by that Gram matrix. In fact, the specification of such an ensemble amounts to merely solving nondegenerate lower triangular linear systems of equations.

We can thus define the shape of a $\varphi$-radial mixture as the coupling of its mixing proportions with the Gram matrix generated by its location vectors:

$$(3.10) \quad [\rho] = (\mathsf{Gram}(\{\tilde{\mu}_k\}), \{q_k\}).$$

We call the two components of this parametrization the *Gram component* and the *mixing component*, respectively. The shape of a profile of $\rho$, say $\breve{\rho}_0(x, y) = \sum_{k=1}^{K} q_k \phi(x, y | HA_0\tilde{\mu}_k)$, corresponding to a rotation $A_0 \in \mathsf{SO}(3)$ will then be given by $[\breve{\rho}_0] = (\mathsf{Gram}(\{HA_0\tilde{\mu}_k\}), \{q_k\})$.

Our interest now is in establishing a relationship between *the shapes of the Radon profiles* of a density and the *shape of the original density*.

**4. Statistical inversion.** Since the coefficients involved in the radial mixture expansion are invariant under projection, the Gram matrix of the location vectors becomes the primary object of interest. Especially in view of sufficiency, we seek a relationship between the Gram components of the projected shape and the original shape. The following theorem provides such a connection and can be seen as an inversion in the shape domain.

THEOREM 4.1 (Shape inversion). *Let $V$ be a $d \times k$ matrix, $2 \leq d < \infty$, $k \leq \infty$, whose columns encode an ensemble of $k$ elements of $\mathbb{R}^d$ with centroid zero. Let $\Psi$ be normalized Haar measure on $\mathsf{SO}(d)$ and $H$ denote the projection onto a subspace of dimension $d - 1$. Then,*

$$(4.1) \qquad \int_{\mathsf{SO}(d)} \mathsf{Gram}\{HAV\}\Psi[dA] = \frac{d-1}{d}\,\mathsf{Gram}\{V\}.$$

PROOF. We may assume that $H = \mathrm{diag}\{1, \ldots, 1, 0\}$ without loss of generality. We notice that $\mathsf{Gram}\{HAV\} = V^\top A^\top HAV$, since $H$ is symmetric idempotent and recognize that $A^\top HA$ is the spectral decomposition of a projection onto the plane $\{A^\top x : x \in \mathrm{Im}(H)\}$, where $\mathrm{Im}(H)$ is the image of $H$. As such, we should be able to encode the same projection relying solely on a unit sphere parametrization, as opposed to using the special orthogonal group. Indeed, $B^\top HB \stackrel{d}{=} I - uu^\top$ for $B \sim \mathrm{Haar}[\mathsf{SO}(d)]$ and $u$ a uniformly random unit vector, $u \sim \mathcal{U}(\mathbb{S}^{d-1})$ ($I - uu^\top$ is the projection onto $u^\perp$). Hence, the proof of the theorem reduces to verifying that, for $u \sim \mathcal{U}(\mathbb{S}^{d-1})$, $\mathbb{E}[uu^\top] = d^{-1}I$. The uniform distribution on the hypersphere is invariant under orthogonal transformations, $Wu \stackrel{d}{=} u$, $\forall W \in \mathsf{O}(d)$. Therefore, $\mathbb{E}[uu^\top] = W\mathbb{E}[uu^\top]W^\top$ for all $W \in \mathsf{O}(d)$, implying that $\mathbb{E}[uu^\top] = cI$ for some constant $c \in \mathbb{R}$. Finally, note that $\mathrm{trace}(\mathbb{E}[uu^\top]) = \mathrm{trace}(\mathbb{E}[u^\top u]) = 1$, so that it must be that $c = d^{-1}$ and the proof is complete. □

The relation in (4.1) reminds one of Cauchy's formula and other related stereological results, where the key ingredient is the isotropy of the projection hyperplanes (see, e.g., Baddeley and Jensen [1]).

Theorem 4.1 says that the expected projected Gram matrix is proportional to the original Gram matrix. Thus, supposing that we can estimate $\{q_k\}$ consistently by some estimator $\{\hat{q}_k\}$, an obvious consistent estimator is given by

$$(4.2) \qquad \left(\frac{d}{(d-1)N} \sum_{n=1}^{N} \mathsf{Gram}(\breve{\rho}_n), \{\hat{q}_k\}\right),$$

which is essentially a method of moments estimator coupled with $\{\hat{q}_k\}$. Unfortunately, things are not so straightforward. Given any profile $\breve{\rho}_n(x, y)$ of $\rho$, the expansion $\breve{p}_n(x,y) = \sum_{k=1}^{K} q_k \phi(x,y | A_n \tilde{\mu}_k)$ is unobservable. Contrary to the case of orthogonal expansions in Hilbert space, there is no transform corresponding to this expansion, so that the unobservable elements $\{q_k\}_{k=1}^{K}$ and $\{A_n \tilde{\mu}_k\}_{1 \le k \le K, 1 \le n \le N}$ must be estimated from the data. A further subtle problem thus arises: since the expansion is unobservable, the correspondence of the indices are also unobservable. When looking at a projection, regardless of how we arrange the location vectors to build the Gram matrix and coefficient vector, the information encoded is the same. However, we must be able to choose this arrangement consistently across all projections, since we will be averaging the Gram matrices across projections. If the indices are unobservable, guaranteeing this consistent construction of the Gram matrices is nontrivial. To surpass this further unidentifiability issue, we impose an assumption on the mixing components.

ASSUMPTION 4.1.  *The components of the density are distinguishable, that is, in the setup given in* (3.6), *we further assume that* $q_i \ne q_j, \ \forall i \ne j$.

Assumption 4.1 allows us to use the auxiliary parameters (the mixing and perhaps the scaling coefficients, if these are included) estimated from the data to recover the unobservable labeling. A hybrid maximum likelihood/method of moments (MoM) estimator is presented in the next section.

4.1. *A hybrid estimator.* In this section, we propose an estimator of the shape of the unknown density when this can be expanded as a finite mixture satisfying Assumption 4.1. The estimator is a hybrid estimator, fusing together a maximum likelihood estimator of the unobservable profile expansions with a method of moments estimator of the final Gram matrix in light of Theorem 4.1. For simplicity and tidiness, we will treat the planar case. The treatment of the $d$-dimensional case, $d \ge 3$, is directly analogous.

In the planar case, we have $\rho : \mathbb{R}^2 \to [0, \infty)$, an essentially bounded density function supported on the disc of radius $\pi$, $\Delta_2 = \{\mathbf{x} \in \mathbb{R}^2 : \|\mathbf{x}\| \le \pi\}$. We let $N$ be a positive integer and $\{A_n\}_{n=1}^{N}$ be independent and identically distributed random elements of the special orthogonal group $\mathsf{SO}(2)$ drawn according to the corresponding normalized Haar measure. Finally, we write $A\rho(\mathbf{x}) := \rho(A^{-1}\mathbf{x})$ for $\mathbf{x} \in \mathbb{R}^2$ and $A \in \mathsf{SO}(2)$. The corresponding stochastic Radon transform is the collection of projections

$$(4.3) \qquad \breve{p}_n(x) := \Pi\{\rho\}(A_n)(x) = \int_{-\infty}^{+\infty} A_n \rho(x,y)\,dy, \qquad x \in [-\pi, \pi].$$

Let $\varphi(\cdot | \tilde{\xi})$ be a planar radial density function centered at $\tilde{\xi}$, and let $\phi(x|0) = \int_{-\infty}^{+\infty} \varphi(x,y|0)\,dy$ be a symmetric one-dimensional location density, centered

at the origin. Our model is

$$(4.4) \qquad \rho(x,y) = \sum_{k=1}^{K} q_k \varphi(x,y|\tilde{\mu}_k), \qquad q_i \neq q_j, \ \forall i \neq j,$$

so that the $n$th profile is $\breve{\rho}_n(x) = \sum_{k=1}^{K} q_k \phi(x|\mu_k^{(n)})$. Here, $\mu_k^{(n)} \in [-\pi, \pi]$ denotes the projection of the $k$th location vector in the $n$th profile of the stochastic Radon transform: $\mu_k^{(n)} = HA_n\tilde{\mu}_k$, $H = (1,0)$. Since we assume that $\rho$ is supported on the disc $\Delta_2$, it must be that diam$\{\text{supp}(\varphi)\} < 2\pi$.

In practice, we observe a discrete version of the profiles, on certain lattice points $\{x_t\}_{t=1}^{T} \in [-\pi, \pi]$, for $T$ a positive integer. In particular, we assume the lattice to be regular, that is, the $x_t$ being equally spaced. Furthermore, the digital images $\{I_n\}_{n=1}^{N}$ of the profiles will be contaminated by additive noise, which is assumed to be Gaussian and white,

$$(4.5) \qquad I_n(x_t) = \breve{\rho}_n(x_t) + \varepsilon_n(t) = \sum_{k=1}^{K} q_k \phi(x_t|\mu_k^{(n)}) + \varepsilon_n(t)$$

for $1 \leq n \leq N$ and $1 \leq t \leq T$. Here, $\{\varepsilon_n(t)\}$ is a collection of $N$ independent Gaussian white noise processes on $\{1, \ldots, T\}$ with variance $\sigma_\varepsilon^2$. By independence, both between and within the white noise processes and the random rotations in (4.3), we may write down the following log-likelihood expression for the parameters of the unobservable mixture expansion:

$$(4.6) \qquad \ell(\mu, q) \propto -\frac{2\pi}{NT} \sum_{n=1}^{N} \sum_{t=1}^{T} \left\| I_n(x_t) - \sum_{k=1}^{K} q_k \phi(x_t|\mu_k^{(n)}) \right\|^2.$$

Maximization of this log-likelihood requires that we choose $NK$ location parameters ($K$ for each profile), as well as a unique set of $K$ mixing proportions to be shared across the profiles, so as to minimize the residual sum of squares between the observed and stipulated profiles (note that $K$ is assumed to be known). Our hybrid estimator for the shape of the two-dimensional density $\rho$ is then formally written as

$$(4.7) \qquad \widehat{[\rho]} = \left( \underbrace{\frac{2}{N} \sum_{n=1}^{N} \mathsf{Gram}(\{\hat{\mu}_k^{(n)}\}_{k=1}^{K})}_{\breve{G}}, \underbrace{\{\hat{q}_k\}_{k=1}^{K}}_{\hat{q}} \right), \qquad (\hat{\mu}, \hat{q}) = \underset{(\mu,q)}{\arg\max}\, \ell(\mu, q).$$

The hybrid estimator is consistent and asymptotically Gaussian. Here, *asymptotic* refers to both important aspects of the problem: the resolution $T$ and the number of profiles $N$. The $T \to \infty$ asymptotics relate to the deconvolution procedure, while the $N \to \infty$ asymptotics relate to the inversion stage. Depending on how one defines the hybrid estimator, there may be an interesting relationship between the two. To state these results, we point out that the underlying probability space is a product space $(\Omega, \mathcal{G}, \mathbb{P}) = (\Omega_1 \times \Omega_2, \sigma(\mathcal{D}_1 \times \mathcal{D}_2), \mathbb{P}_1 \times \mathbb{P}_2)$. The triple $(\Omega_2, \mathcal{D}_2, \mathbb{P}_2)$ is the space of sequences of special orthogonal matrices $\{A_n\}_{n=1}^\infty$ with $\mathbb{P}_2$ such that all finite-dimensional measures are product normalized Haar measures. On the other hand, the space $(\Omega_1, \mathcal{G}_1, \mathbb{P}_1)$ induces infinite sequences of row-wise independent white noise triangular arrays $\{\varepsilon_n(t,T); 1 \leq t \leq T \leq \infty\}_{n=1}^\infty$ under $\mathbb{P}_1$, with common variance $\sigma_\varepsilon^2$.

First, we note that the ML estimators are consistent for every $N \geq 1$ as $T \to \infty$, implying consistency of the mixing component estimator.

THEOREM 4.2 (MLE deconvolution consistency). *Let* $\tilde{\theta}(\omega_2|N)$ *denote the true parameters for the profiles of a stochastic Radon transform of length* $N$ *of the mixture expansion* (4.4),

$$\tilde{\theta}(\omega_2|N) = (\{\mu_k^{(n)}(\omega_2)\}_{1 \leq k \leq K, 1 \leq n \leq N}, \{q_k\}_{1 \leq k \leq K})$$

*and let* $\hat{\theta}(\omega_1, \omega_2|T, N)$ *denote the corresponding maximum likelihood deconvolution estimate based on the observed profile* $\{I_N(x_t)\}_{t=1}^T$. *Then,*

$$(4.8) \qquad \hat{\theta}(\omega_1, \omega_2|T, N) \xrightarrow[T \to \infty]{\mathbb{P}} \tilde{\theta}(\omega_2|N).$$

To prove Theorem 4.2, we first state without proof a variant of a uniform weak law for triangular arrays due to Jennrich [25].

LEMMA 4.1. *Let* $\{X_{t,T}\}_{t \leq T}$ *be a triangular array of random variables with mean zero and variance* $\sigma_X^2 < \infty$, *such that elements of the same row are independent. Let* $\{g_{t,T}(\theta)\}$ *be a triangular array of continuous functions on a compact Euclidean set* $\Theta$ *satisfying*

$$(4.9) \qquad \lim_{T \to \infty} \sup_{\theta_1, \theta_2 \in \Theta} \left| \frac{1}{T} \sum_{t=1}^T g_{t,T}(\theta_1) g_{t,T}(\theta_2) \right| < \infty.$$

*Then,*

$$\sup_{\theta \in \Theta} \left| \frac{1}{T} \sum_{t=1}^T g_{t,T}(\theta) X_{t,T} \right| \xrightarrow[T \to \infty]{P} 0.$$

PROOF OF THEOREM 4.2. It suffices to show that for $\mathbb{P}_2$-almost all $\omega_2$

$$(4.10) \qquad \hat{\theta}(\omega_1, \omega_2|T, N) \xrightarrow[T \to \infty]{\mathbb{P}_1} \tilde{\theta}(\omega_2|N)$$

for then the bounded convergence theorem will imply that $\forall \epsilon > 0$,

$$\lim_{T \to \infty} \mathbb{P}[\|\hat{\theta}(\omega_1, \omega_2|T, N) - \tilde{\theta}(\omega_2|N)\| > \epsilon]$$

$$= \int_{\Omega_2} \lim_{T \to \infty} \mathbb{P}_1[\|\hat{\theta}(\omega_1, \omega_2|T, N) - \tilde{\theta}(\omega_2|N)\| > \epsilon] \, d\mathbb{P}_2 = 0.$$

To this aim, let $\Theta$ be the support of $\theta$, which is by definition compact, and define $Q_T(\cdot)$ as follows:

$$Q_T(\theta) = -\frac{2\pi}{NT} \sum_{n=1}^{N} \sum_{t=1}^{T} \left\| I_n(x_t) - \sum_{k=1}^{K} q_k \phi(x_t - \mu_k^{(n)}) \right\|^2$$

$$= -\frac{2\pi}{NT} \sum_{n=1}^{N} \sum_{t=1}^{T} \| I_n(x_t) - \Phi_n(x_t|\theta) \|^2.$$

If we establish the existence of a deterministic function $Q(\theta)$ such that:

(1) $\sup_{\theta \in \Theta} |Q_T(\theta) - Q(\theta)| \xrightarrow{\mathbb{P}_1} 0$, as $T \to \infty$,
(2) $\sup_{\theta \,:\, \|\theta - \tilde{\theta}\| \geq \epsilon} Q(\theta) < Q(\tilde{\theta})$, $\forall \epsilon > 0$,
(3) $Q_T(\hat{\theta}(T, N)) \geq Q_T(\tilde{\theta}) - o_{\mathbb{P}_1}(1)$,

for $\mathbb{P}_2$-almost all $\omega_2$, then relation (4.10) will immediately follow (e.g., van der Vaart [52]). We fix and suppress $\omega_2$ so that in what follows, random variables are to be seen as functions of $\omega_1$ only. Let

$$(4.11) \qquad Q(\theta) := -\frac{1}{N} \sum_{n=1}^{N} \int_{-\pi}^{\pi} \| \breve{\rho}_n(x) - \Phi_n(x|\theta) \|^2 \, dx - \sigma_\varepsilon^2,$$

where $\breve{\rho}_n(\cdot) = \Phi_n(\cdot|\tilde{\theta})$ is the $n$th profile (without noise contamination) and $\sigma_\varepsilon^2$ is the variance of the noise component. For tidiness write $\delta_n(x|\theta) := |\breve{\rho}_n(x) - \Phi_n(x|\theta)|$ so that $Q(\theta) = -N^{-1} \sum_n \int \delta_n^2(x|\theta) \, dx - \sigma_\varepsilon^2$. Furthermore, let $\delta_n^*(x_t|\theta) := \sup_{x_t - 2\pi/T \leq y \leq x_t} \delta_n(x|\theta)$ and

$$Q_T^*(\theta) := -\frac{2\pi}{NT} \sum_{n=1}^{N} \sum_{t=1}^{T} \delta_n^{*2}(x_t|\theta) - \sigma_\varepsilon^2$$

$$= -\frac{2\pi}{NT} \sum_{n=1}^{N} \sum_{t=1}^{T} \sup_{x_t - 2\pi/T \leq y \leq x_t} \{\delta_n(y|\theta)\}^2 - \sigma_\varepsilon^2,$$

the second equality following form monotonicity of $y \mapsto y^2$ on $\mathbb{R}^+$. To verify condition (1), we must show convergence of $\sup_{\theta \in \Theta} |Q_T(\theta) - Q(\theta)|$ to zero in probability. Using the triangle inequality on the uniform norm, one obtains

$$(4.12) \qquad \sup_{\theta \in \Theta} |Q_T(\theta) - Q(\theta)| \leq \underbrace{\sup_{\theta \in \Theta} |Q_T(\theta) - Q_T^*(\theta)|}_{A(T)} + \underbrace{\sup_{\theta \in \Theta} |Q_T^*(\theta) - Q(\theta)|}_{B(T)}.$$

For the term $B(T)$, we note that $f_T(\theta) = |Q_T^*(\theta) - Q(\theta)|$ is continuous on the compact set $\Theta$. Furthermore, by the definition of the upper Riemann–Stieltjes integral, we have that $f_T \downarrow 0$ pointwise as $T \to \infty$. It follows by Dini's theorem that $f_T$ converges uniformly to zero, that is, $B(T) = \sup_{\theta \in \Theta} |Q_T^*(\theta) - Q(\theta)| \xrightarrow{T \to \infty} 0$. Consider now $A(T)$:

$$A(T) = \frac{2\pi}{NT} \sum_{n=1}^{N} \sum_{t=1}^{T} [(I_n(x_t) - \Phi_n(x_t|\theta))^2 - \delta_n^{*2}(x_t|\theta) - \sigma_\varepsilon^2]$$

$$= \frac{2\pi}{NT} \sum_{n=1}^{N} \sum_{t=1}^{T} [(\breve{\rho}_n(x_t) + \varepsilon_n(x_t) - \Phi_n(x_t|\theta))^2 - \delta_n^{*2}(x_t|\theta) - \sigma_\varepsilon^2]$$

$$= \frac{2\pi}{NT} \sum_{n=1}^{N} \sum_{t=1}^{T} [\delta_n^2(x_t|\theta) - \delta_n^{*2}(x_t|\theta) + (\varepsilon_n^2(x_t) - \sigma_\varepsilon^2)$$
$$+ 2\varepsilon_n(x_t)(\breve{\rho}_n(x_t) - \Phi_n(x_t|\theta))].$$

The triangle inequality yields

$$\sup_{\theta \in \Theta} |A(T)| \leq \sup_{\theta \in \Theta} \left| \frac{2\pi}{NT} \sum_{n=1}^{N} \sum_{t=1}^{T} (\varepsilon_n^2(x_t) - \sigma_\varepsilon^2) \right|$$
$$+ \sup_{\theta \in \Theta} \left| \frac{2\pi}{NT} \sum_{n=1}^{N} \sum_{t=1}^{T} (\delta_n^2(x_t|\theta) - \delta_n^{*2}(x_t|\theta)) \right|$$
$$+ \sup_{\theta \in \Theta} \left| \frac{2\pi}{NT} \sum_{n=1}^{N} \sum_{t=1}^{T} 2\varepsilon_n(x_t)(\breve{\rho}_n(x_t) - \Phi_n(x_t|\theta)) \right|.$$

The first term on the right-hand side is functionally independent of $\theta$ and converges to zero in $\mathbb{P}_1$-probability as $T \to \infty$ by the weak law of large numbers for triangular arrays.

To see that the second term converges to zero, we notice that the function $|\frac{2\pi}{NT} \sum_{n=1}^{N} \sum_{t=1}^{T} (\delta_n^2(x_t|\theta) - \delta_n^{*2}(x_t|\theta))|$ is continuous on the compact set $\Theta$, and, by definition of the upper Riemann–Stieltjes integral,

$$\left| \frac{2\pi}{NT} \sum_{n=1}^{N} \sum_{t=1}^{T} (\delta_n^2(x_t|\theta) - \delta_n^{*2}(x_t|\theta)) \right| \uparrow 0 \qquad \text{as } T \to \infty,$$

pointwise in $\theta$. It therefore follows from Dini's theorem that convergence to zero is uniform.

Hence, it remains to show that the penultimate term converges to zero in $\mathbb{P}_1$-probability. To this aim, we will use Lemma 4.1. To see that it applies here, we need to verify condition (4.9) for the array $\{\psi_{t,T}(\theta)\} = \{(\breve{\rho}(x_t) - \Phi_n(x_t|\theta)\}$. Observe that, by definition of the upper Riemann–Stieltjes integral,

$$\left| \frac{2\pi}{T} \sum_{t=1}^{T} \psi_{t,T}(\theta_1) \psi_{t,T}(\theta_2) \right| \leq \frac{2\pi}{T} \sum_{t=1}^{T} |\psi_{t,T}(\theta_1)| |\psi_{t,T}(\theta_2)|$$

$$\leq \frac{2\pi}{T} \sum_{t=1}^{T} \delta_n^*(x_t|\theta_1)\delta_n^*(x_t|\theta_2) \downarrow \langle \delta_n(x|\theta_1), \delta_n(x|\theta_2)\rangle_2$$

$$< \infty$$

for all $\theta_1, \theta_2 \in \Theta$. Notice that the limit is a bounded function on $\Theta \times \Theta$. Dini's theorem implies that the convergence to this limit occurs uniformly. We have therefore verified that (1) holds.

We now proceed to verify condition (2). When the location parameters $\{\mu_k^{(n)}\}_{k=1}^K$ within a profile are distinct, the mapping $\{(q_k, \mu_k^{(n)})\}_{k=1}^K \mapsto \sum_{k=1}^K q_k \times \phi(\cdot|\mu_k^{(n)})$ is an injection. Since the $\omega_2$ for which the projected means are distinct is a set of $\mathbb{P}_2$-probability 1, condition (2) follows immediately from the fact that $Q(\cdot)$ [defined in (4.11)] is a norm on $L^2[-\pi, \pi]$. Condition (3) is trivially satisfied, by definition of $\hat{\theta}_T$ as the argmax of $Q_T(\cdot)$. Finally, statements in the proof hold for all $\omega_2 \in \Omega_2$, and the proof is complete. $\square$

Using Theorem 4.2, we establish the consistency of the estimator of the Gram component.

THEOREM 4.3 (Gram component consistency). *Let $\rho$ be as in (4.4), and let $\hat{G}(T, N)$ denote the hybrid estimator of the Gram component of $[\rho]$, Gram($[\rho]$), based on $N$ independent profiles $\{I_n(x_t)\}_{t=1}^T$, $n = 1, \ldots, N$. Then, $\hat{G}$ is $L^p$-consistent for Gram($[\rho]$) in the sense that for every $p > 0$, there exists a sequence $T_N \uparrow \infty$, such that*

$$(4.13) \qquad \mathbb{E}\|\hat{G}(T_N, N) - \text{Gram}([\rho])\|_F^p \overset{N\to\infty}{\longrightarrow} 0,$$

*where $\|\cdot\|_F$ is the Frobenius matrix norm.*

Before proceeding with the proof, we give two lemmas without proof.

LEMMA 4.2. *Let $\alpha(T) = (\alpha_1(T), \ldots, \alpha_K(T))$ be a sequence of random measures on $\{1, \ldots, K\}$ and $\pi_T$ be a sequence of random permutations, defined by the property that $\pi_T(\alpha_1(T), \ldots, \alpha_K(T)) = (\alpha_{(1)}(T), \ldots, \alpha_{(K)}(T))$ a.s. If $\alpha = (\alpha_1, \ldots, \alpha_K)$ is a measure with distinct components, then*

$$(4.14) \qquad \alpha_T \overset{P}{\underset{T\to\infty}{\longrightarrow}} \alpha \quad \implies \quad \pi_T \overset{P}{\underset{T\to\infty}{\longrightarrow}} \pi,$$

*where $\pi$ is defined by the property $\pi(\alpha_1, \ldots, \alpha_K) = (\alpha_{(1)}, \ldots, \alpha_{(K)})$.*

LEMMA 4.3. *If $W_T$ is a sequence of random $p \times p$ permutation matrices converging to a permutation matrix $W$ in probability, and $X_T$ is a sequence of random $p \times 1$ vectors converging to $X$ in probability, then $W_T X_T \overset{P}{\underset{T\to\infty}{\longrightarrow}} WX$.*

PROOF OF THEOREM 4.3. Let $\hat{\mu}_k^{(n)}(N, T)$ be the ML estimator of the $k$th location parameter within the $n$th projection. Let $\hat{S}_n(T, N)$ be the corresponding estimate of the Gram matrix for the locations within the $n$th projection, $\hat{S}_n(T, N) = \{\langle \hat{\mu}_i^{(n)}(N, T), \hat{\mu}_j^{(n)}(N, T)\rangle\}_{i,j=1}^K$, and $S_n$ be the true Gram matrix corresponding to the true means in the $n$th projection. Then, we have

$$(4.15)\qquad \hat{G}(T, N) = \frac{2}{N}\sum_{n=1}^N \hat{S}_n(T, N) \quad \text{and} \quad \tilde{G}(N) = \frac{2}{N}\sum_{n=1}^N S_n$$

and we can bound the $L^p$ distance $(\mathbb{E}\|\hat{G}(T, N) - \mathsf{Gram}([\rho])\|_F^p)^{1/p}$ above by

$$\underbrace{(\mathbb{E}\|\hat{G}(T, N) - \tilde{G}(N)\|_F^p)^{1/p}}_{A(T, N)} + \underbrace{(\mathbb{E}\|\tilde{G}(N) - \mathsf{Gram}([\rho])\|_F^p)^{1/p}}_{B(N)}.$$

It is straightforward that $\lim_{N\to\infty} B(N) = 0$, so we concentrate on $A(T, N)$:

$$A(T, N) \le \frac{2}{N}\sum_{n=1}^N (\mathbb{E}\|\hat{S}_n(T, N) - S_n\|_F^p)^{1/p} = 2(\mathbb{E}_2[\mathbb{E}_1\|\hat{S}_1(T, N) - S_1\|_F^p])^{1/p}.$$

Here we have used exchangeability and Fubini's theorem. By Theorem 4.2, we have that the maximum likelihood deconvolution estimates are $\mathbb{P}$-weakly consistent (in the sense of convergence in probability) so that, for every $N \in \mathbb{N}$, we have $(\hat{S}_1(T, N), \ldots, \hat{S}_N(T, N)) \xrightarrow{\mathbb{P}}_{T\to\infty} (S_1, \ldots, S_N)$. This is true by a combination of the continuous mapping theorem for convergence in probability, Lemmas 4.2 and 4.3. It follows by the fact that $\{\hat{S}(T, N)\}_{T\in\mathbb{N}}$ is uniformly bounded so that $\mathbb{E}\|\hat{S}_1(T, N) - S_1\|_F^p \xrightarrow{T\to\infty} 0$, for all $N \in \mathbb{N}$ (e.g., see Corollary 2.2.2, page 38 of Lukacs [36]). Now consider a nonnegative sequence $\{b_n\}$ converging monotonically to zero. Since $\lim_{T\to\infty} A(T, N) = 0$ for any value of $N$, then we can find a sequence $\{T_N\}$ such that $A(T_N, N) \le b_N, \forall N$. Taking the limit as $N \to \infty$ completes the proof. $\square$

Taking the consistency results in Theorems 4.2 and 4.3 as a starting point, we may also prove weak convergence results for the mixing and Gram components.

THEOREM 4.4 (Mixing component CLT). *Let $\hat{q}_{N,T}(\omega_1, \omega_2)$ denote the maximum likelihood estimator based on $N$ independent profiles $\{I_n(x_t)\}_{t=1}^T$, $n = 1, \ldots, N$, of the mixing proportions $q$ of the mixture given in (4.4). Let $N_T \uparrow \infty$ be a sequence dominated by $T$, such that $\hat{q}_{N_T,T}(\omega_1, \omega_2) \xrightarrow{\mathbb{P}} q$ for $T \to \infty$. Then, if $\varphi$ is twice differentiable, it holds that*

$$(4.16)\quad \mathbb{P}[\sqrt{N_T T}(\hat{q}_{N_T,T}(\omega_1, \omega_2) - q) \in (\mathbf{y}, \mathbf{y} + d\mathbf{y})] \xrightarrow{T\to\infty} \frac{\exp\{-\mathbf{y}^\top F\mathbf{y}/2\}}{(2\pi)^{K/2}|F|^{-1/2}}\, d\mathbf{y},$$

*where the entries of the matrix $F$ are given by*

$$(4.17) \qquad F_{ij} = \frac{1}{2\pi\sigma_\varepsilon^2} \mathbb{E}\left\{ \int_{-\pi}^{\pi} \phi(x|\mu_i)\phi(x|\mu_j)\, dx \right\}$$

*with $\phi(x|\mu_j) := \int_{-\infty}^{+\infty} \varphi(x, y| A\tilde{\mu}_j)\, dy$ for $A \sim \mathrm{Haar}[\mathsf{SO}(3)]$, $j = 1, \ldots, K$.*

PROOF. As with the proof of Theorem 4.2, it will suffice to show that for $\mathbb{P}_2$-almost all $\omega_2$, relation (4.16) holds with $\mathbb{P}_1$ replacing $\mathbb{P}$. We start out with a technical note, whose relevance will become clear later. Fix two indices $i, j \leq K$ and consider the collection $\{\int_{-\pi}^{\pi} \phi(y|\mu_i^{(n)})\phi(y|\mu_j^{(n)})\, dy\}_{n \geq 1}$. This comprises an i.i.d. sequence in $L^1(\Omega_2, \mathcal{G}_2, \mathbb{P}_2)$ so that by the strong law of large numbers, the set $B_{i,j}$ defined as

$$\left\{ \omega_2 \in \Omega_2 : \frac{1}{N} \sum_{n=1}^{N} \int_{-\pi}^{\pi} \phi(y|\mu_i^{(n)})\phi(y|\mu_j^{(n)})\, dy \right.$$

$$\left. \xrightarrow{N\to\infty} \mathbb{E} \int_{-\pi}^{\pi} \phi(y|\mu_i^{(n)})\phi(y|\mu_j^{(n)})\, dy \right\}$$

has $\mathbb{P}_2$-probability 1. Hence, $\mathbb{P}[\bigcap_{i,j} B_{i,j}] = 1$. For the rest of the proof, we fix an arbitrary $\omega_2 \in \bigcap_{i,j} B_{i,j}$ which will not be explicitly written out. Theorem 4.2, allows the following Taylor expansion of the gradient of the log-likelihood around the maximum likelihood estimator $\hat{q}_T$ for large $T$:

$$\nabla\ell(\hat{q}_{N_T,T}) = \nabla\ell(q) + [\nabla^2\ell(q)](\hat{q}_{N_T,T} - q) + O_{\mathbb{P}}(\|\hat{q}_{N_T,T} - q\|^2),$$

so that

$$\frac{1}{\sqrt{N_T T}}\nabla\ell(q) + \frac{1}{N_T T}[\nabla^2\ell(q)]\sqrt{N_T T}(\hat{q}_{N_T,T} - q)$$

$$+ \frac{1}{\sqrt{N_T T}}O_{\mathbb{P}}(\|\hat{q}_{N_T,T} - q\|^2) = 0.$$

In order to proceed, we calculate $\nabla\ell$ and $\nabla^2\ell$ as follows:

$$\frac{\partial}{\partial q_i}\ell(q) = -\frac{1}{\sigma_\varepsilon^2}\sum_{n=1}^{N}\sum_{t=1}^{T}\phi(x_t|\mu_i^{(n)})\left(I_n(x_t) - \sum_{k=1}^{K} q_k\phi(x_t|\mu_k^{(n)})\right),$$

$$\frac{\partial^2}{\partial q_i\,\partial q_j}\ell(q) = \frac{1}{\sigma_\varepsilon^2}\sum_{n=1}^{N}\sum_{t=1}^{T}\phi(x_t|\mu_i^{(n)})\phi(x_t|\mu_j^{(n)}).$$

Now, it can be seen that $\nabla^2\ell$ does not depend on $\omega_1 \in \Omega_1$. Fix two indices $i, j \leq K$ and observe that by the triangle inequality,

$$\left| \frac{1}{NT}\sum_{n=1}^{N}\sum_{t=1}^{T}\phi(x_t|\mu_i^{(n)})\phi(x_t|\mu_j^{(n)}) - \frac{1}{2\pi}\mathbb{E}\int_{-\pi}^{\pi}\phi(y|\mu_i)\phi(y|\mu_j)\, dy \right|$$

$$\leq \frac{1}{N} \sum_{n=1}^{N} \underbrace{\left| \frac{1}{T} \sum_{t=1}^{T} \phi(x_t|\mu_i^{(n)})\phi(x_t|\mu_j^{(n)}) - \frac{1}{2\pi} \int_{-\pi}^{\pi} \phi(y|\mu_i^{(n)})\phi(y|\mu_j^{(n)}) \, dy \right|}_{\alpha(n,T)}$$

$$+ \underbrace{\left| \frac{1}{2\pi} \frac{1}{N} \sum_{n=1}^{N} \int_{-\pi}^{\pi} \phi(y|\mu_i^{(n)})\phi(y|\mu_j^{(n)}) \, dy - \frac{1}{2\pi} \mathbb{E} \int_{-\pi}^{\pi} \phi(y|\mu_i)\phi(y|\mu_j) \, dy \right|}_{\beta(N)}.$$

Choose any $\epsilon > 0$. An argument involving Dini's theorem, such as the one in the proof of Theorem 4.2, shows that the first term on the right-hand side converges to zero uniformly in $\mu_i^{(n)}$ and $\mu_j^{(n)}$. Hence, convergence to zero is uniform over $n$, too. As a result, we can choose a $T_0$ such that $\alpha(n,T) < \epsilon/2$ for any $T \geq T_0$ and for all $n \in \mathbb{N}$. Since $\omega_2 \in B$, we can also choose an $N_0$ such that for any $N \geq N_0$ it holds that $\beta(N) < \epsilon/2$. Consequently, for any $\epsilon > 0$ we can choose an $M = T_0 \wedge N_0$ such that

$$\left| \frac{1}{NT} (\nabla^2 \ell)_{ij} - \frac{1}{2\pi\sigma_\varepsilon^2} \mathbb{E} \int_{-\pi}^{\pi} \phi(y|\mu_i)\phi(y|\mu_j) \, dy \right| < \epsilon$$

for $T, N \geq M$. Thus, we have established that, $\mathbb{P}$ almost surely

$$(4.18) \qquad \frac{1}{NT} \{(\nabla^2 \ell)_{ij}\} \xrightarrow{T,N \to \infty} \left\{ \frac{1}{2\pi\sigma_\varepsilon^2} \mathbb{E} \int_{-\pi}^{\pi} \phi(x|\mu_i)\phi(x|\mu_j) \, dx \right\} = \{-F_{ij}\}$$

with $F$ being the Fisher information matrix. This remains true when replacing $N$ by an increasing sequence $N_T$, and take the limit as $T \uparrow \infty$. We now turn to show that the gradient of the log-likelihood satisfies a central limit theorem. Let $N_T \uparrow \infty$ be as in the statement of the theorem. Define a triangular array of random $K$-vectors $\{Y_{T,n}\}$ with $1 \leq n \leq N_T \leq T$

$$(4.19) \qquad Y_{T,n} := \frac{1}{\sqrt{TN_T}} \frac{1}{2\pi\sigma_\varepsilon^2} \sum_{t=1}^{T} \delta(x_{t,T}, n)\vec{\phi}(x_{t,T}, n),$$

where $\{x_{t,T}\}$ is a regular lattice of $T$ points on $[-\pi, \pi]$, and the vectors $\vec{\phi}$ and scalars $\delta$ are defined, respectively, as

$$\vec{\phi}(x_{t,T}, n) := \begin{pmatrix} \phi(x_{t,T}|\mu_1^{(n)}) \\ \vdots \\ \phi(x_{t,T}|\mu_K^{(n)}) \end{pmatrix},$$

$$\delta(x_{t,T}, n) := \frac{(I_n(x_{t,T}) - \sum_{k=1}^{K} q_k \phi(x_{t,T}|\mu_k^{(n)}))}{2\pi\sigma_\varepsilon^2}.$$

We note that $\mathbb{E}_1 Y = 0$ throughout the array. Furthermore, by independence,

$$\text{Cov}_1[Y_{T,n}] = \mathbb{E}_1[Y_{T,n} Y_{T,n}^\top]$$

(4.20)

$$= \frac{1}{TN_T} \frac{1}{2\pi\sigma_\varepsilon^2} \sum_{t=1}^{T} \vec{\phi}(x_{t,T}, n)\vec{\phi}(x_{t,T}, n)^\top < \infty.$$

Since we have fixed $\omega_2 \in \bigcap_{i,j} B_{i,j}$, taking the limit as $T \uparrow \infty$ yields

(4.21)
$$\sum_{n=1}^{N_T} \mathbb{E}_1[Y_{T,n} Y_{T,n}^\top] \xrightarrow{T \to \infty} \frac{1}{2\pi\sigma_\varepsilon^2} \int_{-\pi}^{\pi} \mathbb{E}_2[\vec{\phi}(y)\vec{\phi}(y)^\top]\, dy,$$

where the integral is understood element-wise. We now show that our triangular array satisfies a Lindeberg condition. We have

$$\sum_{n=1}^{N_T} \mathbb{E}_1[\|Y_{T,n}\|^2; \|Y_{T,n}\| > \epsilon]$$

$$\leq \sum_{n=1}^{N_T} \mathbb{E}_1\left[\left\|\frac{1}{\sqrt{TN_T}} \sum_{t=1}^{T} \delta(x_{t,T}, n)c\vec{1}\right\|^2; \left\|\frac{1}{\sqrt{TN_T}} \sum_{t=1}^{T} \delta(x_{t,T}, n)c\vec{1}\right\| > \epsilon\right],$$

where $\vec{1} \in \mathbb{R}^K$ is a vector of 1's and $c > 0$ is an upper bound for all the elements of $\vec{\phi}$. These are differentiable and compactly supported functions of $x \in [-\pi, \pi]$ for every $n$. The $\delta$'s are i.i.d. Gaussian random variables over $n$ and $T$, so, by defining i.i.d. random variables $X_n \overset{d}{=} X \sim \mathcal{N}(0, \sigma_\varepsilon^2 c)$, the last line may be re-written as

$$\frac{1}{N_T} \sum_{n=1}^{N_T} \mathbb{E}_1\left[\|X_n \vec{1}\|^2; \left\|\frac{1}{\sqrt{N_T}} X_n \vec{1}\right\| > \epsilon\right] = \int_{\Omega_1} KX^2 \mathbf{1}_{\{\sqrt{K}|X| > \epsilon\sqrt{N_T}\}}\, d\mathbb{P}_1,$$

the last term converging to 0 by virtue of the dominated convergence theorem and $X$ having finite variance. It follows from the multidimensional Lindeberg–Feller theorem (see, e.g., Proposition 2.27 on page 20 of [52]) that

(4.22)
$$\frac{1}{\sqrt{TN_T}} \nabla \ell = \sum_{n=1}^{N_T} Y_{T,n} \xrightarrow{T \to \infty} \mathcal{N}\left(\mathbf{0}, \frac{1}{2\pi\sigma_\varepsilon^2} \int_{-\pi}^{\pi} \mathbb{E}_2[\vec{\phi}(y)\vec{\phi}(y)^\top]\, dy\right).$$

Returning to the Taylor expansion, as $T \to \infty$ the $(N_T T)^{-1/2}$-scaled error vanishes in probability (by assumption) and Slutsky's lemma yields

(4.23)
$$\sqrt{TN_T}(\hat{q}_{N_T,T} - q) \xrightarrow{T \to \infty} \mathcal{N}_K(\mathbf{0}, \sigma_\varepsilon^2 F^{-1}), \qquad \mathbb{P}_2\text{-a.s.}$$

This completes the proof. □

THEOREM 4.5 (Gram component CLT). *Let $\rho$ be as in (4.4) and let $\hat{G}(T, N)$ denote the estimator of the Gram component of $[\rho]$, $\mathsf{Gram}([\rho])$, based on $N$ independent profiles $\{I_n(x_t)\}_{t=1}^T$, $n = 1, \ldots, N$. Then, there exists* $\tau_N \uparrow \infty$ *such that the normalized difference* $\hat{G}(\tau_N, N) - \mathsf{Gram}([\rho])$ *is asymptotically distributed according to the matrix Gaussian distribution*

$$(4.24) \qquad \sqrt{N}(\hat{G}(\tau_N, N) - \mathsf{Gram}([\rho])) \overset{N \to \infty}{\Longrightarrow} \mathscr{N}_{K \times K}(\mathbf{0}, \Sigma),$$

*where the covariance matrix* $\Sigma$ *is given by*

$$(4.25) \qquad \Sigma = (V^\top \otimes V^\top) \operatorname{Cov}[\operatorname{vec}\{\Gamma\}](V \otimes V).$$

*Here,* $\otimes$ *stands for the Kronecker product,* vec *is the column stacking operator,* $V$ *is any* $3 \times K$ *matrix satisfying* $V^\top V = \mathsf{Gram}([\rho])$ *and* $\Gamma$ *is a* $3 \times 3$ *random matrix with the following second-order properties:*

$$(4.26) \qquad \begin{aligned} \operatorname{var}(\Gamma_{ii}) &= \tfrac{1}{9}, \qquad \operatorname{var}(\Gamma_{ij}) = \tfrac{1}{15} \qquad \textit{for } i \neq j, \\ \operatorname{cov}(\Gamma_{ii}, \Gamma_{jj}) &= -\tfrac{1}{18} \qquad \textit{for } i \neq j \end{aligned}$$

*and with uncorrelated off-diagonal elements.*

We give here the definition of the matrix Gaussian distribution (also see Chapter 2 of Nagar and Gupta [37]).

DEFINITION 4.1. Let $M$ be an $p \times n$ real matrix, and let $\Sigma$ and $\Phi$ be positive definite $p \times p$ and $n \times n$ matrices, respectively. A real random matrix $X$ is said to have the *matrix Gaussian distribution* $\mathscr{N}_{p \times n}(M, \Sigma \otimes \Phi)$ if

$$\operatorname{vec}(X^\top) \sim \mathscr{N}_{pn}(\operatorname{vec}(M^\top), \Sigma \otimes \Phi).$$

PROOF OF THEOREM 4.5. Let $F_{N,T}$ denote the distribution of $\sqrt{N}\hat{G}(T, N)$ and $H$ denote the distribution $\mathscr{N}_{K \times K}(\mathsf{Gram}([\rho]), \Sigma)$, that is, the "shifted" stipulated limiting distribution with expectation $\mathsf{Gram}([\rho])$ instead of zero. If $\tilde{G}(N)$ is as in the proof of Theorem 4.3, we denote by $Q_N$ the distribution of $\sqrt{N}\tilde{G}(N)$. It will suffice to show that, for some $\tau_N \uparrow \infty$, $d_{\mathrm{Pr}}(F_{\tau_N, N}, H) \overset{n \to \infty}{\to} 0$, where $d_{\mathrm{Pr}}$ denotes the Prokhorov metric, metrizing weak convergence (e.g., Billingsley [4]). Applying the triangle inequality,

$$(4.27) \qquad d_{\mathrm{Pr}}(F_{\tau_N, N}, H) \leq d_{\mathrm{Pr}}(F_{\tau_N, N}, Q_N) + d_{\mathrm{Pr}}(Q_N, H).$$

Now let $d_P(X, Y) = \inf\{\delta > 0 : \mathbb{P}[d(X, Y) > \delta] < \delta\}$ be the standard metric that metrizes convergence in probability, for $d$ the metric on the range of the random variables. Now $d_{\mathrm{Pr}}(\Lambda_1, \Lambda_2)$ is the infimum of $d_P(X, Y)$ over all pairs $(X, Y)$ of random variables with $(\Lambda_1, \Lambda_2)$ the respective distributions, provided that $d$ induces a separable space (Dudley [13]). It follows that $d_{\mathrm{Pr}}(F_{\tau_N, N}, Q_N)$ will converge to zero if we can show that

$$(4.28) \qquad \sqrt{N}A(\tau_N, N) = \sqrt{N}\mathbb{E}\|\hat{G}(\tau_N, N) - \tilde{G}(N)\|_F^p \overset{N \to \infty}{\to} 0.$$

For this to hold, it must be that $A(\tau_N, N)$ is $o(N^{-(1+\alpha)/2})$ for some $\alpha > 0$. Since $\lim_{T \to \infty} A(T, N) = 0$ for all $N \in \mathbb{N}$, we can always choose such a $\tau_N$.

We have established that the first term in the right-hand side of (4.27) converges to zero. For the second term, we use the classical central limit theorem. The sequence $\{2 S_n\}$ is i.i.d. with mean $\mathsf{Gram}([\rho])$ and covariance $\Sigma$ (see (4.19) in Panaretos [42], with $n = 2$). Therefore, by the multidimensional central limit theorem, one has $2N^{-1/2} \sum_{n=1}^{N} S_n \overset{N \to \infty}{\Longrightarrow} H$. $\quad \square$

In Theorems 4.3 and 4.5, the order of magnitude of $T$ is made dependent on that of $N$. We discuss this briefly, starting from the setup in Theorem 4.3. When performing the deconvolution, we ask that the mixing proportions are the same for all profiles. As a result, the dimension of the unknown parameter grows with $N$, so that we should make $T$ depend on $N$. If the deconvolutions are performed independently for each profile,

$$(4.29) \quad \underset{(\{q_k^{(n)}\}_{k=1}^K, \{\mu_k^{(n)}\}_{k=1}^K)}{\arg\min} \frac{2\pi}{T} \sum_{t=1}^{T} \left\| I_n(x_t) - \sum_{k=1}^{K} q_k^{(n)} \phi(x_t | \mu_k^{(n)}) \right\|^2,$$
$$n = 1, 2, \ldots, N,$$

then exchangeability implies a "stronger" consistency result

$$(4.30) \quad \forall \epsilon > 0 \ \exists T_0, N_0 : \mathbb{E} \| \hat{G}(T, N) - \mathsf{Gram}([\rho]) \|_F^p < \epsilon$$
$$\forall N \geq N_0 \text{ and } T \geq T_0.$$

A compromise between the two extreme approaches of *overall* and *separate* sums of squares is to assign profiles into groups of size less than a fixed number, not depending on $N$. In the case of overall optimization, more projection data may require better resolution $T$. In biological practice, of course, the instruments will give a certain—hopefully high—resolution. This depends on the current state of technology and can be thought of as being inflexible. On the other hand, the number of projections can become arbitrarily large, the only constraint being computing time (Glaeser [18]).

Theorem 4.5 says that for $T$ converging to infinity sufficiently fast as compared to $N$, we are in a consistent regime, and a classical central limit theorem applies for the estimated shape. A slow growth of $T$, however, may provide asymptotically biased estimators. If deconvolution is carried out separately for each profile, a standard $\sqrt{n}$-type central limit theorem for the maximum likelihood estimates of the projected location parameters applies for the MLE deconvolution within a single profile. A delta method argument subsequently implies that picking $\tau_N$ of the order of $N^\kappa$ for any $\kappa > 1$ is sufficient in order to guarantee weak convergence to the stipulated limit.

4.2. *Exact deconvolution and a theorem of Carathéodory.* Our analysis has underlined the difficulties posed by the unobservability of the radial mixture expansion parameters. Intuitively, there is no transform! It is natural to ask whether there are special cases where these expansions need not be estimated, but are available as a determinate aspect of the data at hand.

Recall that the (deterministic) problem of deconvolution is the solution of a *Fredholm integral equation of the first kind* of the form $\int_{-\pi}^{+\pi} g(t - y)h(y)\,dy = f(t)$, for $g$, when $h$ and $f$ are known. The singular value decomposition of the operator $\phi \mapsto g * \phi$ (e.g., Kanwal [30]) may be used to solve such an equation when the functions involved are square integrable. In the present case, however, the function $g(x) = \sum_{k=1}^{K} \alpha_k \delta(x - \mu_k)$ is a weighted Dirac comb, which is not an element of $L^2$.

Let $f(x)$ be a profile. In the absence of noise, our deconvolution problem is described as $f(x) = \phi(x) * \sum_{k=1}^{K} \alpha_k \delta(x - \mu_k)$, where $f$ and $\phi$ are known, and we wish to recover $\{\mu_k\}_{k=1}^{K}$ and $\{\alpha_k\}_{k=1}^{K}$. In practice, we observe a discretely sampled profile $\{f_t\}_{t=-T/2}^{T/2-1}$, $f_t = f(\frac{2\pi t}{T})$ for $T$ even (say). Therefore, if we apply the inverse discrete Fourier transform we obtain $d_f^{-1}(\kappa)/\hat{\phi}(\kappa) \approx \sum_{k=1}^{K} \alpha_k e^{i\mu_k\kappa}$ for $T$ large, assuming that $\hat{\phi}(\kappa)$ does not vanish for $\kappa \in \{2\pi t/T\}_{t=-T/2}^{T/2}$. This translates our deconvolution problem into one of frequency identification, allowing us to invoke the following theorem.

THEOREM 4.6 (Carathéodory [7]). *Let $\{c_k\}_{k=0}^{n-1}$ be complex constants, at least one of which is nonzero ($n > 1$). Then, there exist an integer $m$, $1 \le m \le n$, real numbers $\beta_j > 0$ and distinct frequencies $\omega_j \in (-\pi, \pi]$ ($j = 1, \ldots, m$) such that the $c_k$ can be uniquely represented as*

$$（4.31) \qquad c_k = \sum_{j=1}^{m} \beta_j e^{i\omega_j k}, \qquad k = 0, \ldots, n-1.$$

Setting $w(\kappa) = d_f^{-1}(\kappa)/\hat{\phi}(\kappa)$, we have that $w(-\kappa) = \overline{w(\kappa)}$, and so we are in the setting $w(\kappa) = \sum_{k=1}^{K} \alpha_k e^{i\mu_k\kappa}$ and want to recover the $\alpha_k$'s and $\mu_k$'s. Carathéodory's theorem assures us that, for $T$ large and for appropriate densities $\phi$, *there exists a unique solution to our deconvolution problem with probability* 1 (owing to our Haar measure assumption). Pisarenko [45] hinged on a constructive proof of Carathéodory's result due to Grenander and Szegö [21] to determine the hidden frequencies $\{\omega_j\}$:

1. Build the matrix $C = \{c_{i-j}\}_{i,j=1}^{n}$.
2. Find the eigenvector $v = (v_0, \ldots, v_m)$ corresponding to the smallest eigenvalue (assuming no multiplicity) of the Toeplitz matrix $C_m$, the top-left submatrix of $C$ of dimension $m \times m$.
3. Find the $K$ distinct unit roots $\{e^{i\omega_j}\}_{j=1}^{m}$ of $\sum_{n=0}^{m} v_n z^n = 0$.

As with any deconvolution problem, this method will be sensitive to the presence of noise. Li and Speed [35], in studying electrophoresis experiments, report that the method is fairly robust to the presence of noise and can be used to provide starting values for a maximum likelihood procedure. Such an approach is unlikely to be succesful in our setting, as, for certain profiles, two (at least) spikes may fall within a critical distance $\epsilon$ of each other, rendering the method very sensitive to noise unless $T$ is overwhelmingly large. The dependence of deconvolution on the relationship between $T$ and $K$ is connected with the *Rayleigh limit* and the problem of *superresolution* (see Donoho et al. [11]). Indeed, we have a *nearly black* object to recover; the $\ell_0$-norm of the unknown signal (spike train) is much smaller than $T$.

The most important aspect in our case, though, is that the information across different profiles can be used to gain insights about the location of the spikes within each particular profile, so that spikes that lie close to one another could potentially be identified. In the noiseless case, Carathéodory's theorem guarantees that we will still recover an expansion, namely that which combines the almost coincident spikes into a single spike.

An extension of the method of Pisarenko to higher dimensions is not straightforward. The one-dimensional result can employed, however, in a coordinate-wise fashion (on the one-dimensional marginals of every projection).

**5. Two examples.** The dimension of the search space and the form of the objective function (4.6) render the practical solution of the optimization problem challenging in its own right. We will not pursue it here, as it is the subject of a separate investigation. However, in order to illustrate both the problem and the application of the hybrid estimator, we consider two mixture data sets, that are purposely chosen to be sparse and noise free, so that Pisarenko's method may be applied.

**5.1.** *A two-dimensional Gaussian mixture.* Assume that we observe a finite sample of $N = 150$ profiles from the stochastic Radon transform of the function

$$\rho(u) = \sum_{j=1}^{5} \frac{j}{2\pi\sigma^2} \exp\left\{ -\frac{1}{2\sigma^2} (u - \mu_j)^\top (u - \mu_j) \right\},$$

(5.1)
$$u = (u_x, u_y) \in \mathbb{R}^2,$$

with $\sigma = 0.3$ and $\mu_1 = (0.6, 0)^\top$, $\mu_2 = (0.6, 0.8)^\top$, $\mu_3 = (-0.1, 0.1)^\top$, $\mu_4 = (-1, -0.3)^\top$, $\mu_5 = (-0.2, -0.6)^\top$. The choice of location parameters and standard deviation was made so as to ensure a certain sparsity. Figure 4(a) gives a contour and intensity plot of the mixture density.

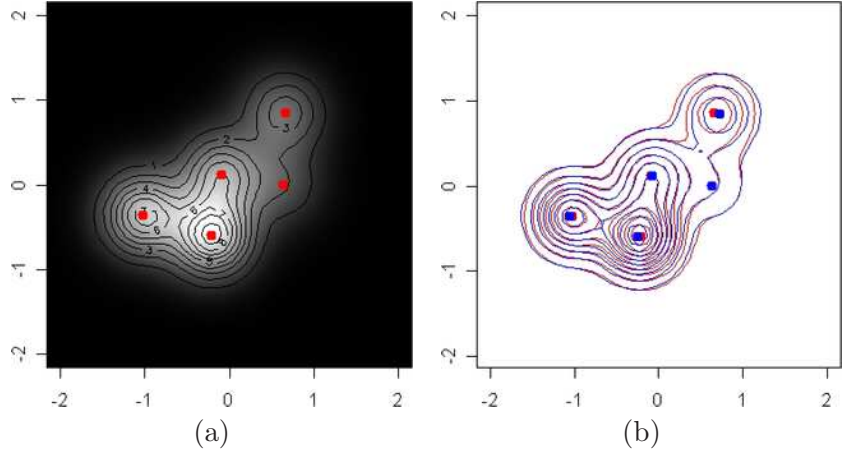

(a)                                      (b)

FIG. 4. (a) *Contour plot of the density ρ superimposed on its intensity plot, with dots indicating the locations of the means.* (b) *Superimposition of the contour plots of the true (red) and estimated (blue) densities.*

The 150 profiles are digitized on a grid of $T = 256$ regular lattice points. A sample of six profiles from this stochastic Radon transform is presented in Figure 5. The projection angles are uknown.

The deconvolution was performed separately for each profile, using only Pisarenko's method, which yielded extremely accurate results, and then the hybrid estimator was constructed. Figure 4(b) contains the contours of the estimated density superimposed on those of the true density.

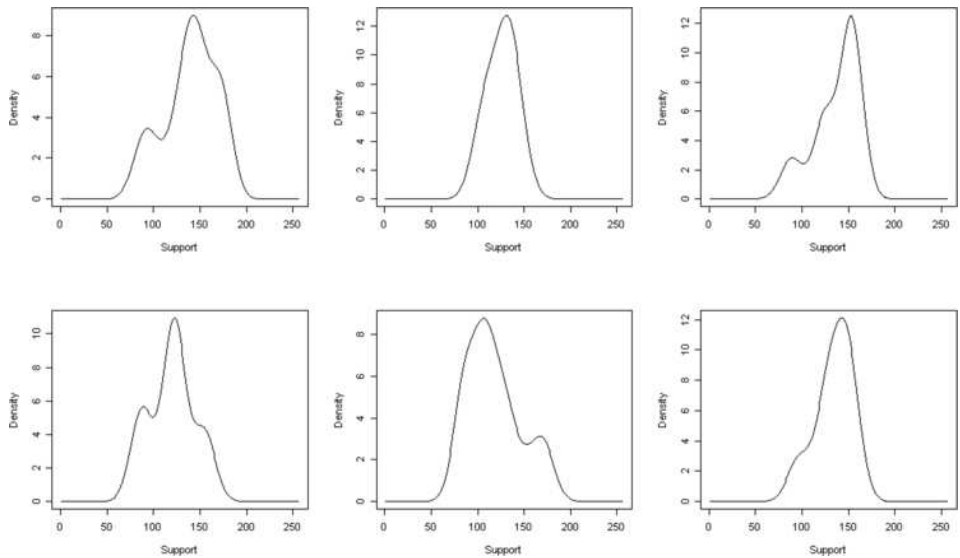

FIG. 5. *Six sample profiles from the 150 profiles of the realisation of the stochastic Radon transform of the density given in* (5.1).

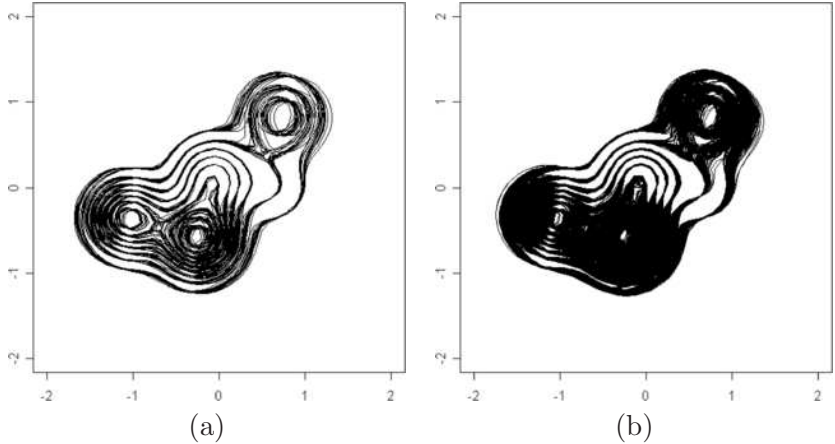

Fig. 6. *Superimposed bootstrap replicates as a means of assessment of uncertainty. Panels* (a) *and* (b) *contain the superimposition of* 15 *and* 100 *replications, respectively.*

A challenging problem is uncertainty estimation and presentation. Motivated by Brillinger et al. [6] and Brillinger, Downing and Glaeser [5], who exploited symmetries to assess variability, we employ a bootstrap approach. We resample the 150 profiles and construct bootstrap replicates of the estimated density. We then superimpose the contour plots (see Figure 6). The rule of thumb is that the more tangled the contours appear the more uncertainty is associated with that particular region. The important aspect of these figures is that the overall shape is seen to be preserved and not to be highly variable.

### 5.2. A three-dimensional Gaussian mixture.

Next, we consider a sparse Gaussian mixture in three dimensions. We observe $N = 150$ profiles from the stochastic Radon transform of the Gaussian mixture

$$\rho(u) = \sum_{k=1}^{4} \frac{q_k}{\sigma^3(\sqrt{2\pi})^3} \exp\left\{ -\frac{(u - \mu_k)^\top (u - \mu_k)}{2\sigma^2} \right\},$$

$$u = (u_x, u_y, u_z) \in \mathbb{R}^3,$$

with $\sigma = 0.46$, $q_1 = 2$, $q_2 = 3$, $q_3 = 2.4$, $q_4 = 4$ and $\{\mu_j\}$ given as $\mu_1 = (0, 0.8, -0.3)^\top$, $\mu_2 = (0.7, -0.4, -0.3)^\top$, $\mu_3 = (-0.7, -0.4, -0.3)^\top$, $\mu_4 = (0, 0, 0.8)^\top$. Visualization of this synthetic particle is challenging since we must visualize level surfaces rather than contours, via an isosurface plot [Figure 8(a)]. Again, we notice that the mixture is sparse in the sense previously described. A sample of four profiles from the realization of the stochastic

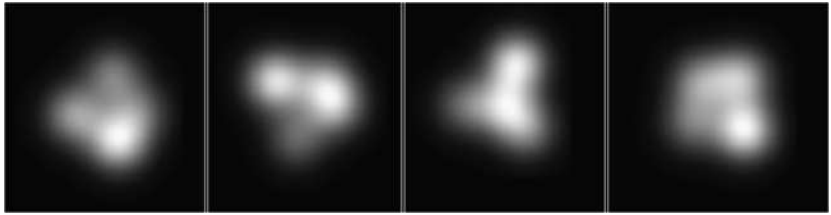

Fig. 7. *Intensity plots for a sample of four profiles from the stochastic Radon transform of ρ. These images would correspond to the data yielded by the electron microscope.*

Radon transform is depicted in Figure 7. The corresponding Euler projection angles are unknown.

Deconvolution was carried out separately for each image, and the mixing coefficients were estimated as before. In particular, we used a naive deconvolution approach, applying Pisarenko's method to the one-dimensional marginals of each projection (which worked well here due to sparsity). A visual comparison of the estimated and the true density is given in Figure 8, where isosurface plots are given for both densities.

**6. Concluding remarks.** We formulated and investigated a problem of statistical tomography where the projection angles are *random* and *unobservable*. The problem was seen to be ill-posed since unobservability of the projection angles limits us to consider inference modulo an appropriate orthogonal group. For essentially bounded and compactly supported densities, these invariants were seen *to be identifiable* and the problem of their recovery was phrased as an estimation problem. The abstraction involved in modu-lar inference may render the pursuit of a solution in the most general case overly ambitious. Indeed it is initially not at all clear how an estimator can be constructed. Nevertheless, it was seen that once the problem is "over-regularized," so as to become an essentially parametric problem, one may develop an appropriate framework for modular inference, and obtain consistent estimators. In particular, assuming that the unknown function admits a radial mixture representation (or, more generally, a radial basis function representation), consistent estimation may be performed hinging on ideas from D. G. Kendall's Euclidean shape theory, without making any attempt to estimate the unknown projection angles. In this setup, a hybrid estimator was constructed whose determination requires a deconvolution step and an inversion step, the former being the more challenging one.

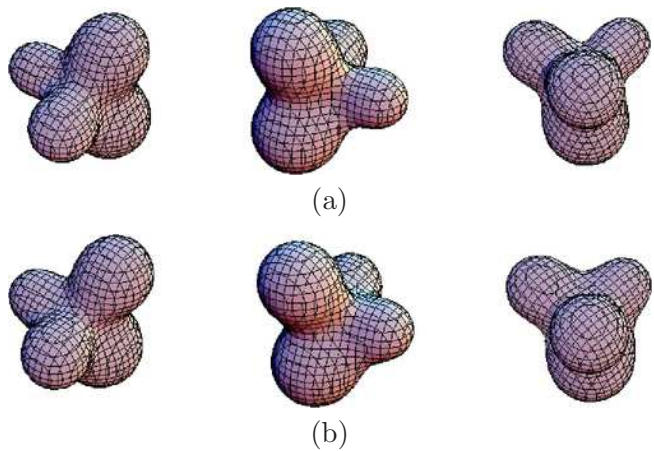

Fig. 8. *(a) Isosurface plots from three different perspective for the density ρ. (b) Isosurface plots from the same three perspecitves of the estimated density.*

Tomography with unknown projection angles arises in singe particle electron microscopy, and biophysicists typically proceed by estimating the unknown angles by means of a prior model (called a *low resolution model*). This low resolution model initializes an iterative procedure that estimates the angles, updates the estimate of the particle and cycles until convergence. Low resolution models often originate from an ad hoc "naked eye examination" of the projections by the experienced scientist, who uses his visual intuition to circumvent the lack of angular information. In fact, many such models comprise an ensemble of solid balls in space, and can be quite successful as starting models, provided that the particle has enough symmetries to enable this ad hoc "naked eye" model construction. The results in this paper suggest that, potentially, it could be practically feasible to construct an "objective" prior model based solely on the data at hand, neither requiring symmetries, nor attempting to estimate the unknown angles. From one point of view, the approach described can be thought of as a mathematical formalization of the biophysicist's visual intuition. Indeed, the radial expansion setup investigated can provide a fruitful framework for the construction of initial models (R. M. Glaeser, personal communication) and its application to the single-particle setup is the subject of ongoing work. On a more theoretical level, the identifiability results presented establish the feasibility of reconstruction from single particle data (up to a coordinate system).

Finally, we add a few comments on the estimation framework. Throughout the paper, the number of mixing components $K$ of the unknown density has been assumed to be known—for example, the scientist will have some insight in its choice. Nevertheless, $K$ could also be treated as an unknown parameter to be estimated, this being a standard problem in mixture estimation (see James, Priebe and Marchette [24], and references therein). In fact, more unknown parameters can be introduced, as long as the mixture remains identifiable. Investigation of how an EM approach could be adapted for this simultaneous deconvolution problem would be of interest in this case.

Another assumption that was made was regarding the isotropy of the Eulerian projection angles. In the anisotropic case, one can hardly proceed at all when the projection angles are unobservable; any re-weighting scheme would be ill-defined. The determination of an algorithm that would perform the simultaneous deconvolution step in a real setting is a problem that is of interest in its own right. Connections with discrete sparse inverse problems, such as those studied in Donoho et al. [11], are especially relevant here.

**Acknowledgments.** I wish to thank Professor David Brillinger for motivating me to work on this problem and for many stimulating discussions. I am indebted to Professor Robert Glaeser, LBNL, for our most helpful interactions on the structural biology side and for his warm hospitality during my visits at the MRC Laboratory of Molecular Biology at Cambridge and the Lawrence Berkeley National Lab. I also wish to thank Professor Sir David Cox for a number of thought-provoking conversations and Professor Anthony Davison for useful comments. Finally, my thanks go to an Associate Editor and the referees for their thorough and constructive comments.

Research primarily carried out while the author was at the Department of Statistics, University of California, Berkeley.

Part of this research was carried out while the author was visiting at the Department of Statistics, Stanford University and the Lawrence Berkeley National Laboratory. Their hospitality is gratefully acknowledged.

## REFERENCES

[1] BADDELEY, A. and VEDEL-JENSEN, E. B. (2005). *Stereology for Statisticians.* Chapman and Hall/CRC, Boca Raton, FL. MR2107000
[2] BERAN, R., FEUERVERGER, A. and HALL, P. (1996). On nonparametric estimation of intercept and slope distributions in random coefficient regression. *Ann. Statist.* **24** 2569–2592. MR1425969
[3] BOOKSTEIN, F. L. (1978). *The Measurement of Biological Shape and Shape Change. Lecture Notes in Biomathematics* **24**. Springer, New York.
[4] BILLINGSLEY, P. (1968). *Convergence of Probability Measures.* Wiley, New York. MR0233396
[5] BRILLINGER, D. R., DOWNING, K. H. and GLAESER, R. M. (1990). Some statistical aspects of low-dose electron imaging of crystals. *J. Statist. Plann. Inference* **25** 235–259. MR1064428
[6] BRILLINGER, D. R., DOWNING, K. H., GLAESER, R. M. and PERKINS, G. (1989). Combining noisy images of small crystalline domains in high resolution electron microscopy. *J. Appl. Statist.* **16** 165–175.
[7] CARATHÉODORY, C. and FEJÉR, L. (1911). Über den zusammenhang der extemen von harmonischen funktionen mit ihren koeffizienten und über den Picard–Landausch Sätz. *Rend. Circ. Mat. Palermo* **32** 218–239.
[8] CHANG, I.-S. and HSIUNG, C. A. (1994). Asymptotic consistency of the maximum likelihood estimate in positron emission tomography and applications. *Ann. Statist.* **22** 1871–1883. MR1329172

[9] CHIU, W. (1993). What does electron cryomicroscopy provide that X-ray crystallography and NMR spectroscopy cannot? *Ann. Rev. Biophys. Biomol. Struct.* **22** 233–255.

[10] DEANS, S. R. (1993). *The Radon Transform and Some of Its Applications.* Krieger, Malabar, FL. MR1274701

[11] DONOHO, D. L., JOHNSTONE, I. M., HOCH, J. C. and STERN, A. S. (1992). Maximum entropy and the nearly black object (with discussion). *J. Roy. Statist. Soc. Ser. B* **54** 41–81. MR1157714

[12] DRENTH, J. (1999). *Principles of Protein X-Ray Crystallography.* Springer, New York.

[13] DUDLEY, R. M. (1968). Distances of probability measures and random variables. *Ann. Math. Statist.* **39** 1563–1572. MR0230338

[14] FEUERVERGER, A. and VARDI, Y. (2000). Positron emission tomography and random coefficients regression. *Ann. Inst. Statist. Math.* **52** 123–138. MR1771484

[15] FRANK, J. (1999). *Three-Dimensional Electron Microscopy of Macromolecular Assemblies.* Academic Press, San Diego.

[16] GARCZAREK, F., DONG, M., TYPKE, D., WITKOWSKA, E., HAZEN, T. C., NOGALES, E., BIGGIN, M. D. and GLAESER, R. M. (2007). Octomeric pyruvate-ferredoxin oxidoreducatse from Desulfovibrio vulgaris. *J. Struct. Biol.* **159** 9–18.

[17] GLAESER, R. M. (1985). Electron crystallography of biological macromolecules. *Ann. Rev. Phys. Chem.* **36** 243–275.

[18] GLAESER, R. M. (1999). Review: Electron crystallography: Present excitement, a nod to the past, anticipating the future. *J. Struct. Biol.* **128** 3–14.

[19] GLAESER, R. M., CHIU, W., FRANK, J., DEROSIER D., BAUMEISTER, W. and DOWNING, K. (2007). *Electron Crystallography of Biological Macromolecules.* Oxford Univ. Press, New York.

[20] GREEN, P. J. (1990). Bayesian reconstructions from emission tomography data using a modified EM algorithm. *IEEE Trans. Med. Imaging* **9** 84–93.

[21] GRENANDER, U. and SZEGÖ, G. (1958). *Toeplitz Forms and Their Applications.* Univ. California Press, Berkeley. MR0094840

[22] HELGASON, S. (1980). *The Radon Transform.* Birkhäuser, Boston. MR0573446

[23] HENDERSON, R. (2004). Realizing the potential of electron cryo-microscopy. *Q. Rev. Biophys.* **37** 3–13.

[24] JAMES, L. F., PRIEBE, C. E. and MARCHETTE, D. J. (2001). Consistent estimation of mixture complexity. *Ann. Statist.* **29** 1281–1296. MR1873331

[25] JENNRICH, R. I. (1969). Asymptotic properties of nonlinear least squares estimators. *Ann. Math. Statist.* **40** 633–643. MR0238419

[26] JENSEN, S. R. (2004). Sufficient conditions for the inversion formula for the $k$-plane Radon transform in $\mathbb{R}^n$. *Math. Scand.* **94** 207–226. MR2053741

[27] JOHNSTONE, I. M. and SILVERMAN, B. W. (1990). Speed of estimation in positron emission tomography and related inverse problems. *Ann. Statist.* **18** 251–280. MR1041393

[28] JONES, M. C. and SILVERMAN, B. W. (1989). An orthogonal series density estimation approach to reconstructing positron emission tomography images. *J. Appl. Statist.* **16** 177–191.

[29] KALLENBERG, O. (2002). *Foundations of Modern Probability.* Springer, New York. MR1876169

[30] KANWAL, R. P. (1997). *Linear Integral Equations: Theory and Technique.* Birkhäuser, Boston. MR1427946

[31] KENDALL, D. G. (1977). The diffusion of shape. *Adv. in Appl. Probab.* **9** 428–430.

[32] KENDALL, D. G. and KENDALL, W. S. (1980). Alignments in two-dimensional random sets of points. *Adv. in Appl. Probab.* **12** 380–424. MR0569434

[33] KENDALL, W. S. and LE, H. (2009). Statistical Shape Theory. In *New Perspectives in Stochastic Geometry* (W. S. Kendall and I. S. Molchanov, eds.). Oxford Univ. Press (forthcoming).

[34] LE, H. and KENDALL, D. G. (1993). The Riemannian structure of Euclidean shape spaces: A novel environment for statistics. *Ann. Statist.* **21** 1225–1271. MR1241266

[35] LI, L. and SPEED, T. (2000). Parametric deconvolution of positive spike trains. *Ann. Statist.* **28** 1270–1301. MR1805784

[36] LUKACS, E. (1975). *Stochastic Convergence.* Academic Press, New York. MR0375405

[37] NAGAR, D. K. and GUPTA, A. K. (2000). *Matrix Variate Distributions.* Chapman and Hall/CRC, Boca Raton, FL. MR1738933

[38] NATTERER, F. (2001). *The Mathematics of Computerized Tomography. Society for Industrial and Applied Mathematics (SIAM)* **32**. Philadelphia, PA. MR1847845

[39] O'SULLIVAN, F. (1995). A study of least squares and maximum likelihood for image reconstruction in positron emission tomography. *Ann. Statist.* **23** 1267–1300. MR1353506

[40] O'SULLIVAN, F. and PAWITAN, Y. (1993). Multidimensional density estimation by tomography. *J. Roy. Statist. Soc. Ser. B* **55** 509–521. MR1224413

[41] PANARETOS, V. M. (2006). The diffusion of Radon shape. *Adv. in Appl. Probab.* **38** 320–335. MR2264947

[42] PANARETOS, V. M. (2008). Representation of Radon shape diffusions via hyperspherical Brownian motion. *Math. Proc. Cambridge Philos. Soc.* **145** 457–470. MR2442137

[43] PANARETOS, V. M. (2008). On random tomography in structural biology. Technical Report, No. 2008-3, Dept. Statistics, Stanford Univ.

[44] PETER, A. and RANGARAJAN, A. (2006). Shape analysis using the Fisher–Rao Riemannian metric: Unifying shape representation and deformation. In *3rd IEEE International Symposium on Biomedical Imaging: Macro to Nano* **1–3** 1164–1167. Arlington.

[45] PISARENKO, V. F. (1973). The retrieval of harmonics from a covariance function. *Geophys. J. R. Astr. S.* **33** 347–366.

[46] SCHERVISH, M. J. (1995). *Theory of Statistics.* Springer, New York. MR1354146

[47] SHEPP, L. A. and KRUSKAL, J. B. (1978). Computerized tomography: The new medical X-ray technology. *Amer. Math. Monthly* **85** 420–439. MR1538734

[48] SHEPP, L. A. and VARDI, Y. (1982). Maximum likelihood reconstruction in positron emission tomography. *IEEE Trans. Med. Imaging* **1** 113–122.

[49] SILVERMAN, B. W., JONES, M. C., WILSON, J. D. and NYCHKA, D. W. (1990). A smoothed EM approach to indirect estimation problems with particular reference to stereology and emission tomography (with discussion). *J. Roy. Statist. Soc. Ser. B* **52** 271–324. MR1064419

[50] SRIDECHADILOK, B., FRASER, C. S., HALL, R. J., DOUDNA, J. A. and NOGALES, E. (2005). Structural roles for human translation initiation factos eIF3 in initiation of protein synthesis. *Science* **310** 1513–1515.

[51] SMALL, C. G. and LE, H. (2002). The statistical analysis of dynamic curves and sections. *Pattern Recogn.* **35** 1597–1609.

[52] VAN DER VAART, A. W. (1998). *Asymptotic Statistics.* Cambridge Univ. Press, Cambridge. MR1652247

[53] VARDI, Y., SHEPP, L. A. and KAUFMAN, L. (1985). A statistical model for positron emission tomography (with discussion). *J. Amer. Statist. Assoc.* **80** 8–37. MR0786595

[54] YOUNES, L. (1998). Computable elastic distances between shapes. *SIAM J. Appl. Math.* **58** 565–586. MR1617630

INSTITUT DE MATHÉMATIQUES
ECOLE POLYTECHNIQUE FÉDÉRALE
DE LAUSANNE
SWITZERLAND
E-MAIL: victor.panaretos@epfl.ch